\documentclass[11pt,a4paper]{article}

\usepackage[T1]{fontenc}
\usepackage{lmodern}
\usepackage[margin=26mm]{geometry}
\usepackage{amsmath,amssymb,amsthm,mathtools}
\usepackage{enumitem}
\usepackage{microtype}
\usepackage{tikz}
\usetikzlibrary{arrows.meta,positioning,calc}
\usepackage{caption}
\usepackage{booktabs,tabularx,array,longtable}
\usepackage{flafter}
\usepackage{needspace}
\usepackage{placeins}
\usepackage{aliascnt}
\usepackage[hidelinks]{hyperref}
\usepackage[capitalise,nameinlink,noabbrev]{cleveref}
\hypersetup{
 pdftitle={Gallai's path decomposition conjecture for K5-minor-free graphs},
 pdfauthor={Junhao Zhang, Qinghai Liu, Yanmei Hong},
 pdfsubject={Gallai path decompositions},
 pdfkeywords={Gallai conjecture, path decomposition, graph minor, planar graph}
}

\allowdisplaybreaks[2]
\setlist[enumerate]{topsep=4pt,itemsep=3pt,leftmargin=2.3em}
\newtheorem{theorem}{Theorem}[section]
\newaliascnt{lemma}{theorem}
\newtheorem{lemma}[lemma]{Lemma}
\aliascntresetthe{lemma}
\newaliascnt{proposition}{theorem}
\newtheorem{proposition}[proposition]{Proposition}
\aliascntresetthe{proposition}
\newaliascnt{corollary}{theorem}
\newtheorem{corollary}[corollary]{Corollary}
\aliascntresetthe{corollary}
\theoremstyle{definition}
\newaliascnt{definition}{theorem}

\aliascntresetthe{definition}
\crefname{theorem}{Theorem}{Theorems}
\crefname{lemma}{Lemma}{Lemmas}
\crefname{proposition}{Proposition}{Propositions}
\crefname{corollary}{Corollary}{Corollaries}
\Crefname{theorem}{Theorem}{Theorems}
\Crefname{lemma}{Lemma}{Lemmas}
\Crefname{proposition}{Proposition}{Propositions}
\Crefname{corollary}{Corollary}{Corollaries}

\DeclareMathOperator{\pn}{pn}
\newcommand{\floor}[1]{\left\lfloor #1\right\rfloor}
\newcommand{\ceil}[1]{\left\lceil #1\right\rceil}

\tikzset{
 vertex/.style={circle,fill=black,draw=black,inner sep=1.55pt},
 newvertex/.style={circle,fill=white,draw=black,line width=.65pt,inner sep=1.65pt},
 edge/.style={line width=.65pt},
 pale/.style={draw=black!25,line width=.5pt},
 selected/.style={draw=black,line width=1.25pt},
 lab/.style={font=\small},
 arrow/.style={-{Stealth[length=2.2mm]},line width=.7pt}
}

\newcommand{\Exc}{\mathcal K}
\newcommand{\pat}[1]{C_{\mathrm{#1}}}
\newcommand{\pV}{C_V}
\newcommand{\pVp}{C_V^{\prime}}
\newcommand{\pU}{C_U}
\newcolumntype{Y}{>{\raggedright\arraybackslash}X}
\definecolor{pathblue}{RGB}{20,72,130}
\definecolor{pathred}{RGB}{158,50,45}
\tikzset{oldpath/.style={draw=pathblue,line width=1.3pt},newpath/.style={draw=pathred,line width=1.3pt,densely dashed}}

\title{\Large\bfseries Gallai's path decomposition conjecture\\[2pt]for $K_5$-minor-free graphs}
\author{Junhao Zhang$^{a,}$\thanks{Corresponding author.
E-mail: 9820240101@nankai.edu.cn (J.~Zhang),
qliu@fzu.edu.cn (Q.~Liu), yhong@fzu.edu.cn (Y.~Hong).},
Qinghai Liu$^{b,c}$, Yanmei Hong$^d$}
\date{}

\begin{document}
\maketitle
\vspace{-1.6em}
\begin{center}
\begin{minipage}{.97\textwidth}
\centering\small\itshape
$^a$School of Mathematical Sciences and LPMC, Nankai University,\\
Tianjin 300071, China\\[2pt]
$^b$Fujian Science \& Technology Innovation Laboratory for\\
Optoelectronic Information of China\\[2pt]
$^c$Center for Applied Mathematics of Fujian Province, China\\[2pt]
$^d$School of Mathematics and Statistics, Fuzhou University,\\
Fuzhou, Fujian 350108, China
\end{minipage}
\end{center}
\vspace{.3em}
\begin{abstract}
Gallai conjectured that the edges of every connected graph on $n$
vertices can be partitioned into at most $\ceil{n/2}$ paths.
We prove that every connected $K_5$-minor-free graph other than $K_3$
and $K_5-e$ admits such a partition into at most $\floor{n/2}$ paths.
The two exceptions require $\ceil{n/2}$ paths. In a smallest
counterexample we find a $4$-connected planar part containing four
vertices of degree five. Their neighborhoods meet the boundary in at
most one vertex. This allows us to apply the planar configurations
of Blanch\'e, Bonamy, and Bonichon and to realize the reduced graph
as a minor of the original graph.
\end{abstract}
\noindent\textbf{Keywords:} path decomposition; Gallai's conjecture;
$K_5$-minor-free graph; planar graph; subdivision.

\section{Introduction}\label{sec:intro}
All graphs are finite and simple. We use the notation of
Bondy and Murty~\cite{BondyMurty}.
A \emph{path decomposition} of a graph $G$ is a set of nontrivial paths
whose edge sets partition $E(G)$. The least number of these paths is the
\emph{path number}, denoted by $\pn(G)$.
Gallai conjectured that every connected graph $G$ satisfies
\[
 \pn(G)\le \ceil{|V(G)|/2}.
\]
The conjecture appears in Lov\'asz~\cite{Lovasz}.
We call $G$ a \emph{Gallai graph} if $\pn(G)\le\floor{|V(G)|/2}$.

We prove the following bound.
\begin{theorem}\label{thm:main}
Let $G$ be a connected $K_5$-minor-free graph.
If $G$ is not isomorphic to $K_3$ or $K_5-e$, then
\begin{equation}\label{eq:main}
 \pn(G)\le\floor{|V(G)|/2}.
\end{equation}
\end{theorem}
The two exceptions are necessary. A triangle needs two paths.
The graph $K_5-e$ has nine edges, whereas each of its paths has at most
four edges. Thus it needs three paths, and three paths suffice.

An \emph{odd semi-clique} is obtained from $K_{2k+1}$ by deleting at most
$k-1$ edges, where $k\ge1$. It has more than $2k^2$ edges.
Since a path has at most $2k$ edges, such a graph is not Gallai.
Bonamy and Perrett~\cite{BonamyPerrett} asked whether these are the only
connected graphs that are not Gallai.

For planar graphs, Blanch\'e, Bonamy, and Bonichon~\cite{BBB} proved that
the only connected graphs that fail the floor bound are $K_3$ and $K_5-e$.
The same exceptions occur for connected graphs of treewidth at most
three~\cite{TreewidthThree} and, more generally, for connected
$3$-degenerate graphs~\cite{ZhangLiuHong}. Liu, Zhang, and
Hong~\cite{LiuZhangHong} studied the corresponding problem for
$K_{3,3}$-minor-free graphs. By Wagner's theorem~\cite{Wagner}, a graph is planar precisely when
it has neither a $K_5$ minor nor a $K_{3,3}$ minor. \Cref{thm:main}
allows a $K_{3,3}$ minor while retaining the planar bound.

We first show that a smallest counterexample has at most one vertex
of degree at most four (\Cref{sec:low}). We then complete a set of at
most three boundary vertices to a clique to obtain a $4$-connected
planar part. In this part we find four vertices of degree five, three
with no boundary neighbor and the fourth with at most one
(\Cref{sec:structure}).

The planar configurations of~\cite{BBB} can now be used at these four
vertices. The choice of the vertices lets us replace the only possible
added edge in the configuration by a path in the original graph
(\Cref{sec:transfer}). Deleting the four vertices and applying the
local reductions gives a simple minor $H$. A decomposition of $H$
into paths and cycles can be lifted at the cost of two new paths
(\Cref{sec:restore}). To obtain the floor bound, we rule out triangle
components of $H$ and show that every $K_5-e$ component has an edge
that is subdivided when the deleted vertices are restored
(\Cref{sec:exceptional,sec:proof}).

\subsection{Notation}
For $Y\subseteq V(G)$, let $N_G(Y)$ be the set of vertices outside $Y$
with a neighbor in $Y$. Write
\[
 V_i(G)=\{v:d_G(v)=i\},\qquad V_{\le i}(G)=\{v:d_G(v)\le i\}.
\]
Deleting an edge set leaves the vertex set unchanged.
The notation $G+\binom{S}{2}$ means that all missing edges within $S$
are added. For vertices $a,b$ on a path $Q$, let $aQb$ denote the
corresponding subpath.

To \emph{suppress} a vertex of degree two is to delete it and add an edge
between its neighbors. We also suppress chains of such vertices.
We use this operation only when the result is simple.
The new edge is called an \emph{auxiliary edge}.
To expand it is to replace it by the original path.

A subdivision replaces each edge by a path. These paths have pairwise
disjoint interiors. The original vertices are called \emph{roots}, and
the replacement paths are called \emph{branch paths}.
Let $M_4$ be the multigraph obtained from a $4$-cycle by doubling two
opposite edges. Let $W_4$ be the wheel with a rim of length four.
Their subdivisions will always be subgraphs of a simple graph.
A \emph{mixed decomposition} is an edge partition into simple paths and
simple cycles. Isolated vertices need no members.

Set $\Exc=\{K_3,K_5-e\}$. These are the only $K_5$-minor-free odd
semi-cliques. Indeed, for an odd semi-clique of order $2k+1$ with
$k\ge3$, choose one endpoint of each deleted edge. At most $k-1$
vertices are chosen. The remaining vertices include a clique of order
at least $k+2\ge5$. The cases $k=1,2$ give exactly the members of
$\Exc$ after $K_5$ is excluded.
In~\cite{ZhangLiuHong}, odd semi-cliques are called $K$-graphs.
Thus its exceptional $K$-components are precisely the members of
$\Exc$ in our setting.

Throughout the proof, a \emph{smallest counterexample} means a connected
$K_5$-minor-free graph $G\notin\Exc$ that violates \eqref{eq:main}
and has the fewest vertices. Write $n=|V(G)|$.
Every smaller connected $K_5$-minor-free graph outside $\Exc$ is then
Gallai.

\section{Vertices of degree at most four}\label{sec:low}
We prove the following lemma.
\begin{lemma}\label{lem:low}
A smallest counterexample $G$ satisfies $|V_{\le4}(G)|\le1$.
\end{lemma}

\subsection{A cut vertex}
\begin{lemma}\label{lem:cut}
There are no subgraphs $G_1,G_2$ such that
\[
 G=G_1\cup G_2,\qquad V(G_1)\cap V(G_2)=\{z\},\qquad
 |V(G_i)|\ge3.
\]
\end{lemma}
\begin{proof}
Both sides are connected. Put $n_i=|V(G_i)|$, so $n_1+n_2=n+1$.
First suppose that neither side belongs to $\Exc$.
By minimality, each $G_i$ has at most $\floor{n_i/2}$ paths.
Their total is at most $\floor{n/2}$ unless both $n_i$ are even.

In that case, add a leaf at $z$ in each $G_i$.
Each resulting graph is smaller than $G$, is $K_5$-minor-free, and is
not exceptional. It has a decomposition into at most $n_i/2$ paths.
Delete the leaf edge. Either a whole path is removed, or a path now
ends at $z$. If a whole path is removed on either side, the count
suffices. Otherwise join the two paths ending at $z$.
The sides meet only at $z$, so the joined path is simple.
Again the total is at most $\floor{n/2}$.

Each exceptional graph has an optimal decomposition with at least two
paths ending at any prescribed vertex $z$. The following lists show this:
\[
\begin{array}{c|lll}
 K_3&zab&zb&\\
 K_5-zw&zabcw&zbwac&zc\\
 K_5-pq&az&abpzq&paqbz.
\end{array}
\]
If just one side is exceptional, choose a path through $z$ on the other
side. If $z$ is internal, split this path at $z$ and join the two pieces
to two endpoint paths of the exceptional side. If $z$ is an endpoint,
join one pair. In either case, one path is saved.
If both sides are exceptional, join two pairs and save two paths.
The total is at most $\floor{n/2}$ in both cases, a contradiction.
\end{proof}

\subsection{Deleting a path and suppressing vertices}
Suppose that $R$ is obtained by deleting edges of $G$, and
$S\subseteq V_2(R)$ with $|S|\le2$.
We call the suppression of $S$ \emph{valid} if it gives a simple graph.
Thus each auxiliary edge has distinct ends. It duplicates neither a
retained edge nor another auxiliary edge.
If the two vertices of $S$ are adjacent, they are suppressed together
along a path of length three. We then delete isolated vertices.

We use two results of Zhang, Liu, and
Hong~\cite[Lemmas 2.5--2.6]{ZhangLiuHong}.
Let a path and a cycle be edge-disjoint and meet in at least one and at
most five vertices. Their union has a decomposition into two paths,
with one exception: $K_5-e$ with a possibly empty pendant path at each
end of the missing edge.
Also, two edge-disjoint cycles whose intersection has between one and
four vertices can be decomposed into two paths.

\begin{lemma}\label{lem:delete}
In a smallest counterexample, there is no path or cycle $F$ and set
$S\subseteq V(F)\cap V_2(G-E(F))$ with the following properties.
If $F$ is a cycle, its length is three or four.
We have $|S|\le2$, and valid suppression of $S$, followed by deletion of
isolated vertices, gives a graph $H$ with at most $n-2$ vertices.
\end{lemma}
\begin{proof}
The graph $H$ is a minor of $G$, so it is $K_5$-minor-free.
Put $n'=|V(H)|\le n-2$.
Expanding its auxiliary edges recovers the nontrivial part of
$G-E(F)$. Different edges insert different vertices.
Thus expansion preserves simple paths.

Every component outside $\Exc$ has a floor-bound decomposition by
minimality. An exceptional $K_5-e$ component with an auxiliary edge
expands to a proper subdivision, which needs only two paths.
Here is a direct check for the two edge types of $K_5-ab$, with the
other vertices named $c,d,e$.
If $cd$ is replaced by $czd$, use
\[
 adbecz,\qquad bcaedz.
\]
If $ac$ is replaced by $azc$, use
\[
 zadceb,\qquad aedbcz.
\]
Further subdivisions do not increase the number of paths.

A $K_5-e$ component without an auxiliary edge has a decomposition into
a Hamilton cycle and a path of length four.
Once a covered path meets this component, combine that path with the
cycle and keep the four-edge path. This costs at most two extra paths.
The path--cycle exception cannot occur here. Among the five vertices,
the old path can use at most the missing edge of $K_5-e$.
The exception would require four chords of the cycle.
Processing other components does not add edges on these five vertices.

An expanded triangle component is a cycle of length three, four, or
five. A cycle of length three or four can be absorbed into a path that
meets it, at a cost of one extra path. A cycle of length five uses both
suppressed vertices. Hence there is at most one such cycle.

First suppose that $F$ is a path. If an expanded $5$-cycle $Z$ occurs,
process it first. It meets $F$. Unless $F\cup Z$ is the path--cycle
exception, replace their union by two paths.
Suppose that the exception occurs, with core $K=K_5-pq$.
The cycle $Z$ is a component after the deleted path is removed and the
auxiliary edges are expanded. Hence every edge from $K$ to the rest
of $G$ lies on $F$ and is incident with $p$ or $q$.
No component of $G-V(K)$ meets both $p$ and $q$.
Such a component would give a $p$--$q$ path outside $K$, and hence a
$K_5$ subdivision. By \Cref{lem:cut}, at most one outside vertex can
be attached at either end. Thus $G$ is $K_5-e$ with at most one leaf
at each end of its missing edge. Without leaves it is exceptional.
With leaves it has six or seven vertices and a three-path
decomposition. Indeed, extend endpoint paths at the two odd-degree
vertices of $K_5-e$. None of these graphs is a counterexample.

The other exceptional components can now be absorbed as above.
Every remaining component meets $F$, since $G$ is connected.
The paths already formed still cover every edge of $F$.
Thus a path meeting the next component is always available.

Now suppose that $F$ is a cycle of length three or four.
If $G-E(F)$ has no edges, then $G=F$ and the claim follows directly.
Otherwise choose a remaining part that meets $F$.
If an expanded $5$-cycle occurs, choose it first.
If the chosen part has a path decomposition, combine one of its
paths with $F$. If it is an expanded triangle, use the two-cycle
result. If it is an unexpanded $K_5-e$, combine its Hamilton cycle
with $F$ and retain its four-edge path.
In each case the relevant cycles meet in at most four vertices.
All later expanded triangles have length at most four, so the
remaining components can be absorbed as before.

The cost is one for $F$, one for each triangle component, and two
for each $K_5-e$ component. Consequently,
\[
 \pn(G)\le1+\sum_{J\text{ a component of }H}\floor{|V(J)|/2}
 \le1+\floor{n'/2}\le\floor{n/2},
\]
a contradiction.
\end{proof}

We also use the following form of the restoration lemma
in~\cite[Lemma 2.8]{ZhangLiuHong}.
\begin{lemma}\label{lem:spare}
Let $E_0=\{u_iw_i:1\le i\le t\}$ be a matching, where $1\le t\le2$
and $d_G(u_i)=3$ for every $i$.
Let $F\subseteq G-E_0$ be a path, or a cycle of length at most four,
that passes internally through every $u_i$.
Put $R=G-E_0-E(F)$ and choose
\[
 S\subseteq V(F)\cap V_2(R),\qquad |S|\le1.
\]
Validly suppress $S$ in $R$ and delete isolated vertices.
Let the resulting graph be $H$.
If $|V(H)|\le n-2$ and every $w_i$ has odd degree in $H$, then this
configuration cannot occur in a smallest counterexample.
\end{lemma}
The cited lemma does not require the original graph to be
$3$-degenerate. Its conditions apply here: $G$ is not an odd
semi-clique, and every non-Gallai component of the smaller minor $H$
is a member of $\Exc$. The lemma also covers endpoints $w_i$ in such
components. Since $|S|\le1$, valid suppression is equivalent to
$R[N_R(S)\cup S]$ being acyclic, as required in the cited lemma.
For the parity check, we use
\begin{equation}\label{eq:spareparity}
 d_H(w_i)\equiv d_G(w_i)-d_F(w_i)-1\pmod2.
\end{equation}
Suppression preserves the degrees of retained vertices.
Thus $d_G(w_i)$ should be even if $w_i$ is outside $F$ or internal on
$F$, and odd if it is an endpoint of $F$.

\subsection{A vertex of degree at most two}
\begin{lemma}\label{lem:two}
If a smallest counterexample has a vertex of degree at most two,
then every other vertex has degree at least five.
\end{lemma}
\begin{proof}
Suppose first that $u,v$ both have degree at most two.
Take a shortest $u$--$v$ path. At each endpoint of degree two, add its
other incident edge. The result is a path or a cycle of length three
or four. It contains every edge incident with $u$ or $v$.
Apply \Cref{lem:delete} with $S=\varnothing$.

Now let $d(u)\le2$ and $d(v)\in\{3,4\}$.
Take a shortest path $Q=uQv$.
If $d(u)=2$, let $t$ be the unused neighbor of $u$ and extend $Q$ by
$tu$. If $d(u)=1$, make no extension at $u$.
All unused endpoint neighbors lie outside $Q$.

Suppose that $d(v)=3$, with unused neighbors $a,b$.
If $ab\notin E(G)$, suppress $v$ after deleting the extended path.
The vertex $u$ becomes isolated, so \Cref{lem:delete} applies.
If $ab\in E(G)$ and $t\in\{a,b\}$, relabel so that $t=b$ and delete
$abuQv$. The remaining neighbors of $v$ are $a,b$, and the edge $ab$
has been deleted. Again apply \Cref{lem:delete}.

It remains that $ab\in E(G)$ and $t$ is absent or outside $\{a,b\}$.
Temporarily remove $vb$. Extend the path at $v$ by $va$, and also by
$ab$ if $d_G(b)$ is odd. The path passes internally through $v$.
Equation~\eqref{eq:spareparity} makes $b$ odd in the reduced graph.
The vertices $u,v$ disappear, so \Cref{lem:spare} applies.

Suppose that $d(v)=4$, with unused neighbors $a,b,c$.
If these include a nonedge $ab$, extend at $v$ by $vc$.
If the two extensions meet, then $t=c$ and $Q$ has length at most two.
The result is a $3$- or $4$-cycle. Otherwise it is a path.
After deletion, $u$ is isolated and $v$ can be suppressed onto $ab$.

If $a,b,c$ form a clique and $t$ is absent or outside this clique,
extend at $v$ by $vabc$. If $t=c$, delete $bcuQva$.
The other choices of $t$ are symmetric.
The remaining neighbors of $v$ are $b,c$, and $bc$ has been deleted.
Again \Cref{lem:delete} gives a contradiction.
\end{proof}

\subsection{Two vertices of degree three}
One of the reductions at two vertices of degree three requires the
following edge addition.
\begin{lemma}\label{lem:addedge}
Let $J$ be $K_5$-minor-free. Let $a,b,c,d$ be distinct vertices with
\[
 ac,ad,bc,bd,cd\in E(J),\qquad ab\notin E(J).
\]
If $J-\{c,d\}$ has no $a$--$b$ path, then $J+ab$ is $K_5$-minor-free.
\end{lemma}
\begin{proof}
We use the following observation. If two $K_5$-minor-free graphs meet
in a clique of size at most three and have no other common vertices,
their union is $K_5$-minor-free. Indeed, in a model of $K_5$ at least
two branch sets avoid the clique. All branch sets avoiding it must
lie on the same side, since they are pairwise adjacent. Restrict the
other branch sets to that side. Each restriction stays connected
because its vertices in the clique are pairwise adjacent, and every
pair of branch sets is still adjacent. This would give a $K_5$ minor
on that side.

Let $A$ and $B$ be the components of $J-\{c,d\}$ containing $a$ and
$b$, respectively. Start with a $K_4$ on $\{a,b,c,d\}$. Attach
$J[A\cup\{c,d\}]$ along the triangle $acd$ and $J[B\cup\{c,d\}]$
along $bcd$. Each remaining component, together with $c,d$, attaches
along $cd$. The resulting graph is $J+ab$, so the observation applies.
\end{proof}

The next lemma restores the six edges at the deleted vertices.
\begin{lemma}\label{lem:sixedges}
Let $ab\in E(H)$, and let $a,b,c,d$ be distinct vertices of $H$.
Obtain $G$ by deleting $ab$ and adding two new vertices $u,v$ and the
six edges $au,uc,ud,bv,vc,vd$. Then $\pn(G)\le\pn(H)+1$.
\end{lemma}
\begin{proof}
Let $P=rPabPs$ be the path containing $ab$ in a decomposition of $H$.
If one of $c,d$, say $c$, is outside $P$, replace $P$ by
\[
 rPaucvbPs,\qquad udv.
\]
Otherwise, reversal and relabeling leave two possible orders:
\[
\begin{array}{c|ll}
 \text{order on }P&\multicolumn{2}{c}{\text{replacement paths}}\\
 r\cdots c\cdots d\cdots ab\cdots s&rPcuaPdv&udPcvbPs\\
 r\cdots c\cdots ab\cdots d\cdots s&uaPcvdPs&rPcudPbv.
\end{array}
\]
Each new path uses disjoint sections of $P$ and distinct new vertices.
The two paths cover $E(P)\setminus\{ab\}$ and the six new edges exactly
once. All other paths remain unchanged.
\end{proof}

\begin{figure}[tbp]
\centering
\begin{tikzpicture}[scale=.9]
\begin{scope}
 \coordinate(a)at(-2,0);\coordinate(u)at(-.8,0);\coordinate(v)at(.8,0);
 \coordinate(b)at(2,0);\coordinate(c)at(0,1.25);\coordinate(d)at(0,-1.25);
 \draw[oldpath](a)--(u)--(c)--(v)--(b);
 \draw[newpath](u)--(d)--(v);
 \foreach \x/\pos in {a/left,b/right,u/below left,v/below right,c/above,d/below}
  \node[vertex,label=\pos:$\x$]at(\x){};
 \node[lab]at(0,-2){six edges in $G$};
\end{scope}
\draw[arrow](2.8,0)--(4,0)node[midway,above,lab]{$-u,-v$};
\begin{scope}[xshift=6cm]
 \draw[oldpath](-1.3,0)--(1.3,0);
 \node[vertex,label=left:$a$]at(-1.3,0){};
 \node[vertex,label=right:$b$]at(1.3,0){};
 \node[vertex,label=above:$c$]at(0,1.25){};
 \node[vertex,label=below:$d$]at(0,-1.25){};
 \node[lab]at(0,-2){auxiliary edge $ab$ in $H$};
\end{scope}
\end{tikzpicture}
\caption{The edge replacement in \Cref{lem:sixedges}. When $c\notin V(P)$,
expand $ab$ along the solid path and add the dashed path $udv$.
Other edges are omitted.}
\label{fig:six}
\end{figure}
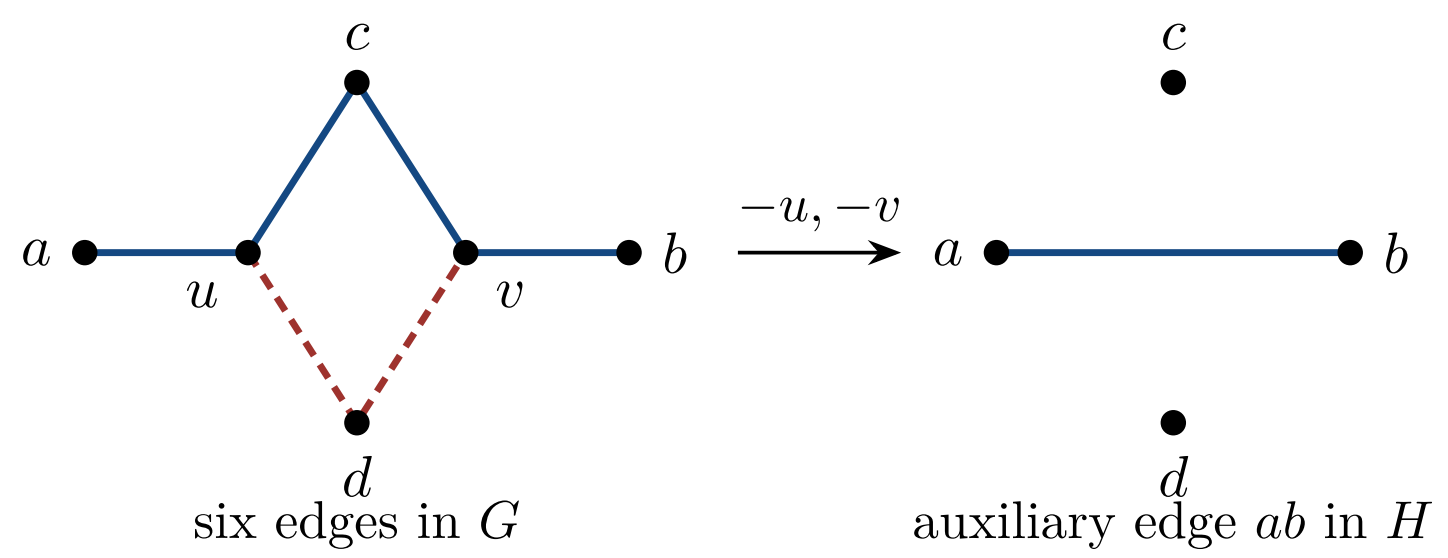

\begin{lemma}\label{lem:three}
A smallest counterexample has no two vertices of degree three.
\end{lemma}
\begin{proof}
Let $d(u)=d(v)=3$, and take a shortest path $Q=uQv$.
Let $A,B$ be the two unused neighbor sets at $u,v$.
These sets lie outside $Q$. If they meet, then $Q$ has length at most two.

\emph{Case 1: $A\cap B=\varnothing$.}
At either endpoint $x$, let $a,b$ be its unused neighbors.
If $ab$ is absent, leave $x$ as an endpoint and later suppress it.
If $ab$ is present, temporarily remove $xb$ and extend by $xa$.
Also use $ab$ if $d_G(b)$ is odd.
The two extensions are disjoint, so they form one simple path.
If both roots are suppressed, use \Cref{lem:delete}.
Otherwise the removed edges form a matching, at most one root is
suppressed, and \eqref{eq:spareparity} gives the required odd degrees.
Use \Cref{lem:spare}.

\emph{Case 2: $A=B=\{a,b\}$.}
Suppose first that $ab$ is absent.
If $d_G(b)$ is even, remove $ub$, take $F=auQv$, and suppress $v$.
If only $a$ has even degree, interchange $a,b$.
If both degrees are odd, remove $ub,va$ and take $F=auQvb$.
Each choice satisfies \Cref{lem:spare}.

Now let $ab$ be present. If $Q=uv$, delete the path $ubav$.
The remaining chain $a-u-v-b$ can be suppressed, since $ab$ was deleted.
If $Q=uwv$, then $N(u)=N(v)=\{a,b,w\}$.
If this set is not a clique, choose $Q$ through its third vertex so
that the unused pair is a nonedge. This is the preceding case.
If it is a clique, delete $ubawv$.
The remaining pairs are $\{a,w\}$ and $\{a,b\}$, and both corresponding
edges were deleted. Apply \Cref{lem:delete}.

\emph{Case 3: $A\cap B=\{c\}$.}
Write $A=\{a,c\}$ and $B=\{b,c\}$, where $a\ne b$.
The following table gives all but one choice. Here $E_0$ is removed
before $F$, and $S$ is the set to suppress.
\begin{center}\small
\begin{tabularx}{\textwidth}{Ylll}
\toprule
Condition&$F$&$E_0$&$S$\\\midrule
$ac,bc$ absent&$uQv$&$\varnothing$&$\{u,v\}$\\
$ac$ absent, $bc$ present; $d(c)$ even&$uQvb$&$\{vc\}$&$\{u\}$\\
Same edges; $d(c)$ odd&$uQvbc$&$\{vc\}$&$\{u\}$\\
$ac,bc$ present; $d(c)$ even&$bcauQv$&$\{uc\}$&$\{v\}$\\
$ac,bc$ present; $d(c)$ odd, $d(b)$ even&$auQvc$&$\{uc,vb\}$&$\varnothing$\\
$a,b,c$ odd; $ac,bc,ab$ present&$cbauQv$&$\{uc\}$&$\{v\}$\\
\bottomrule
\end{tabularx}\end{center}
Interchange $u,v$ and $a,b$ for the symmetric cases.
Each suppressed pair is a nonedge or an edge deleted by $F$.
When $E_0\ne\varnothing$, equation~\eqref{eq:spareparity} gives odd
degrees at the required vertices. Thus \Cref{lem:delete,lem:spare}
apply. The only remaining case is
\begin{equation}\label{eq:lastthree}
 ac,bc\in E(G),\qquad ab\notin E(G),\qquad d(a),d(b),d(c)\text{ odd}.
\end{equation}
If $Q=uv$, put $H=G-\{u,v\}+ab$.
Contracting $au,bv$ gives $H$ after parallel edges are removed.
The graph is connected, since $a,c,b$ are connected and contain all
outside neighbors of $u,v$.
The vertices $a,b,c$ are odd in $H$, so $H\notin\Exc$.
Use a floor-bound decomposition of $H$. Replace $ab$ in its path by
$auvb$ and add $ucv$. This uses at most one extra path.

Hence $Q=udv$. The vertices $c,d$ are common neighbors of $u,v$.
Apply the same case analysis to the path $ucv$.
We may therefore also assume $ad,bd\in E(G)$ and $d_G(d)$ odd.
If $cd$ is absent, remove $vb$, take $F=uacvdb$, and suppress $u$.
Its remaining pair is $\{c,d\}$, while $b$ is an odd-degree endpoint.
Use \Cref{lem:spare}.

Finally, suppose that $cd$ is present.
If $G-\{u,v,c,d\}$ contains an $a$--$b$ path $T$, delete
$F=uaTbdcv$. The remaining pairs at $u,v$ are $\{c,d\}$ and $\{b,d\}$.
Their edges lie on $F$, so \Cref{lem:delete} applies.
Otherwise put $J=G-\{u,v\}$ and $H=J+ab$.
By \Cref{lem:addedge}, $H$ is $K_5$-minor-free.
The four neighbors are connected in $J$ and contain all outside
neighbors of $u,v$. Thus $H$ is connected.
The vertices $a,b,c,d$ are odd in $H$, so $H\notin\Exc$.
Minimality and \Cref{lem:sixedges} give
\[
 \pn(G)\le\pn(H)+1\le\floor{(n-2)/2}+1=\floor{n/2},
\]
a contradiction.
\end{proof}

\subsection{Degrees three and four}
\begin{lemma}\label{lem:threefour}
A smallest counterexample has no pair of vertices of degrees three
and four.
\end{lemma}
\begin{proof}
Let $d(u)=3$ and $d(v)=4$. Choose a shortest path $Q=uQv$.
Its unused neighbor sets $A,B$ have sizes two and three.
They lie outside $Q$.
We use the following extension at a degree-four endpoint.
If its three unused neighbors contain a nonedge $bc$, extend through
the third neighbor $a$. If they form a clique, extend by $vabc$.
In both cases $v$ is internal and its remaining pair $\{b,c\}$ can
be validly suppressed.

If $A\cap B=\varnothing$, use this extension at $v$.
If $A$ is a nonedge, leave $u$ as an endpoint and suppress both roots.
If $A=\{a,b\}$ with $ab\in E(G)$, remove $ub$ and extend at $u$ by
$ua$, followed by $ab$ when $d_G(b)$ is odd.
Apply \Cref{lem:spare}, suppressing only $v$.

Suppose that $A\cap B=\{c\}$, with $A=\{a,c\}$ and $B=\{b,d,c\}$.
If $ac$ is absent, use the degree-four extension at $v$ and suppress
both roots. Their remaining pairs differ since $a\notin B$.
Now let $ac$ be present. If $bd$ is absent, delete $uQvca$.
If $bd$ is present and one of $bc,dc$ is present, rename so that $bc$
is present and delete $uQvdbca$.
The remaining pairs are, respectively, $\{a,c\},\{b,d\}$ and
$\{a,c\},\{b,c\}$. Every required edge has been deleted.
Use \Cref{lem:delete}.
If $bd$ is present but $bc,dc$ are absent, remove $uc$.
Take $F=auQvb$ when $d_G(c)$ is even, and $F=cauQvb$ when it is odd.
The vertex $u$ disappears, $v$ has the nonedge pair $\{c,d\}$,
and $c$ is odd in the reduced graph. Use \Cref{lem:spare}.

It remains that $A\subseteq B$. Write $A=\{a,b\}$ and $B=\{a,b,c\}$.
If $ac$ is absent, delete $uQvb$ when $ab$ is absent, and $uQvba$
when $ab$ is present. The case where $bc$ is absent is symmetric.
Otherwise $ac,bc$ are present. Delete $uQvbca$ if $ab$ is absent,
and $uQvbac$ if $ab$ is present.
In each case the remaining pairs are distinct, and their edges are
absent or deleted. Apply \Cref{lem:delete}.
\end{proof}

\subsection{Two vertices of degree four}
\begin{lemma}\label{lem:four}
A smallest counterexample has no two vertices of degree four.
\end{lemma}
\begin{proof}
Let $d(u)=d(v)=4$, and choose a shortest path $Q=uQv$.
Let $A,B$ be its unused neighbor sets, each of size three.
We delete one path or a cycle of length three or four, and then
suppress $u,v$.
If $A\cap B=\varnothing$, use the degree-four extension from the
preceding proof at both ends. The extensions are disjoint and the
remaining pairs are distinct.

Suppose that $A\cap B=\{c\}$, with $A=\{a,b,c\}$ and $B=\{d,e,c\}$.
If neither set is a clique, choose a nonedge in each.
Let $x,y$ be the unused third vertices and delete $xuQvy$.
If $x=y$, this is a cycle of length three or four, since $Q$ has
length at most two. The two suppressed pairs differ.
Otherwise assume that $B$ is a clique. Start with $F=auQvdec$ and
append $cb$ if $bc$ is present.
The remaining pairs are $\{b,c\}$ and $\{e,c\}$.
Their edges are absent or lie on $F$.

Suppose that $A\cap B=\{b,c\}$, with $A=\{a,b,c\}$ and $B=\{d,b,c\}$.
If neither $a$ nor $d$ sees both $b,c$, choose a nonedge incident with
$a$ and one incident with $d$. Extend through the respective third
vertices $x,y$, and delete $xuQvy$.
If $x=y$, it is again a short cycle.
The auxiliary edges differ because they contain $a$ and $d$.
Otherwise assume $db,dc\in E(G)$.
Take $F=dbuQvc$ and append $ca$ if $ac$ is present.
The remaining pairs $\{a,c\}$ and $\{d,b\}$ are distinct and valid.

Finally, let $A=B=\{a,b,c\}$.
If at least two edges are absent, rename so that $ab,ac$ are absent
and delete $cuQvb$.
If exactly one edge is absent, call it $ac$ and delete $cuQvba$.
The remaining pairs are $\{a,b\}$ and $\{a,c\}$ in both cases.
If $a,b,c$ form a clique, contract the interior of $Q$ towards $u$,
keeping $v$. This gives a $K_5$ on $\{u,v,a,b,c\}$, a contradiction.

Every listed suppression is valid and removes at least two vertices.
Each is excluded by \Cref{lem:delete}.
\end{proof}

\begin{proof}[Proof of \Cref{lem:low}]
A pair with one degree at most two is excluded by \Cref{lem:two}.
The remaining degree pairs are $(3,3)$, $(3,4)$, and $(4,4)$.
They are excluded by \Cref{lem:three,lem:threefour,lem:four}.
\end{proof}

\section{Four vertices of degree five in a planar part}\label{sec:structure}

The planar reductions take place near four vertices of degree five.
We first find a planar subgraph after completing a set of at most
three boundary vertices to a clique. We then choose the four vertices
so that their neighborhoods contain at most one boundary vertex.

\subsection{Choosing the planar part}
\begin{lemma}\label{lem:part}
Let $G$ be a connected $K_5$-minor-free graph with at least two vertices.
Suppose that $d_G(x)\ge5$ for every $x$ other than a specified vertex $z$.
There is a nonempty connected set $X\subseteq V(G)\setminus\{z\}$ such
that, for $S=N_G(X)$, we have $1\le |S|\le3$ and
\[
 P=G[X\cup S]+\binom{S}{2}
\]
is $4$-connected and planar. Each missing edge $ab$ within $S$ has an
$a$--$b$ path in $G$ whose internal vertices lie outside $X\cup S$.
\end{lemma}
\begin{proof}
Consider all nonempty connected sets $D\subseteq V(G)\setminus\{z\}$
with $1\le |B|\le3$, where $B=N_G(D)$, for which
$J=G[D\cup B]+\binom{B}{2}$ is $K_5$-minor-free and every missing
edge within $B$ has a path as in the statement. Any component of $G-z$
is such a set, with boundary $\{z\}$. Choose $D$ of minimum order.
Every vertex of $D$ has the same neighbors in $G$ and $J$, so
$|V(J)|\ge6$.

Suppose that $T$ is a minimum vertex cut of $J$ and $|T|\le3$.
The clique $B\setminus T$ meets at most one component of $J-T$.
Choose a component $C$ disjoint from $B$, and write
$J_C=G[C\cup T]+\binom{T}{2}$. Every component of $J-T$ has a
neighbor at every vertex of $T$, by the minimality of the cut.
In particular, $N_G(C)=T$, since all edges incident with $C$ are
original. Also $T$ meets $D$: deleting only vertices of $B$ leaves
$D$ connected, and every remaining boundary vertex has a neighbor
in $D$. Thus $C$ is a proper subset of $D$.

For each missing edge $ab$ of $G[T]$, another component of $J-T$
gives an $a$--$b$ path in $J$ whose interior avoids $C\cup T$.
Replace each added edge of $B$ on this path by its prescribed path
in $G$. The new internal vertices lie outside $D\cup B$, so the
resulting walk still avoids $C\cup T$ internally. It contains the
required $a$--$b$ path. Hence, whenever $J_C$ is $K_5$-minor-free,
$C$ contradicts the choice of $D$.

If $|T|=1$, then $J_C$ is a subgraph of $J$. If $|T|=2$, contract
a path through another component between the two vertices of $T$ to
an edge and discard the remaining outside vertices. This gives $J_C$
as a minor of $J$. Both cases contradict minimality, so $J$ is
$3$-connected.

Suppose that $J$ is nonplanar. It has a vertex of degree at least five
and hence differs from the cubic Wagner graph. By the small-cut form
of Wagner's theorem due to Ma, Steffen, Wolf, and
Zhang~\cite[Lemma~2.4]{MaSteffenWolfZhang}, $J$ has a three-vertex cut
$T=\{a,b,c\}$ for which $J-T$ has at least three components.
Choose $C$ as above. Contract a second component into $a$ to obtain
$ab,ac$, and a third into $b$ to obtain $bc$. This shows that $J_C$
is a minor of $J$, again a contradiction.

Thus $J$ is planar. If it has a three-vertex cut $T$, contract another
component of $J-T$ to a single vertex with neighborhood $T$ and
delete all other vertices outside $C\cup T$ and this vertex. In a
plane drawing, replace the three edges at the contracted vertex by
a triangle on $T$, drawing its sides along the three edges and keeping
one copy of each side. This gives a plane drawing of $J_C$, so $J_C$
is $K_5$-minor-free. This last contradiction proves that $J$ is
$4$-connected. Take $X=D$, $S=B$, and $P=J$.
\end{proof}

\subsection{Choosing the four vertices}
\begin{theorem}\label{thm:structure}
Let $G$ be a connected $K_5$-minor-free graph with at least two vertices
and $|V_{\le4}(G)|\le1$. There is a set $X$ as in \Cref{lem:part}
containing four distinct degree-five vertices $u_1,u_2,u_3,u_4$ such that
\begin{equation}\label{eq:roots}
 N_G(u_i)\cap S=\varnothing\quad(1\le i\le3),
 \qquad |N_G(u_4)\cap S|\le1.
\end{equation}
\end{theorem}
\begin{proof}
Choose the possible vertex of degree at most four as $z$, or choose
$z$ arbitrarily if all degrees are at least five. Take $X,S,P$ from
\Cref{lem:part} and fix a plane embedding of $P$. Each $x\in X$ has
the same neighbors in $P$ and $G$. Every triangle of $P$ bounds a face,
since a triangle with vertices on both sides would be a three-vertex
cut. Moreover, $P$ contains no $K_4$: any further vertex lies in
one face of a drawing of the $K_4$, and the vertices on that face
separate it from the fourth vertex.

Write $s=|S|$, $b=|E_P(X,S)|$, and let $t$ be the number of degree-five
vertices in $X$. Put $\varepsilon=3|V(P)|-6-|E(P)|$. Euler's formula
gives $\varepsilon=\sum_f(|\partial f|-3)\ge0$. Since $S$ is a clique,
$\sum_{v\in S}d_P(v)=b+s(s-1)$. Consequently
\[
 t\ \ge\ \sum_{x\in X}(6-d_P(x))
   = b+(s-3)(s-4)+2\varepsilon.
\]
Let $r$ count the vertices of $X$ with two neighbors in $S$.
A vertex with three neighbors in $S$ would give a $K_4$.
Hence precisely $b-r$ vertices of $X$ have a neighbor in $S$.
If $t_0$ counts the degree-five vertices in $X$ with no such neighbor,
then
\[
 t_0\ge t-(b-r)\ge (s-3)(s-4)+r+2\varepsilon.
\]

For $s=1$ this gives $t_0\ge6$. For $s=2$, the $r$ common neighbors
of the boundary pair correspond to the triangular faces at its edge.
The remaining $2-r$ incident faces each contribute at least one to
$\varepsilon$. Thus $r\le2$ and
$t_0\ge 2+r+2(2-r)\ge4$. In both cases take four vertices counted by $t_0$.

Now let $s=3$. The triangle $S$ bounds a face. The three faces on the
other sides of its edges are distinct. Indeed, two such faces could
coincide only if the two corresponding boundary edges were consecutive
on both sides at their common vertex, which has degree at least four.
Exactly $r$ of these faces are triangular: no vertex outside $S$ sees
all three of its vertices. Therefore $r\le3$ and
$\varepsilon\ge3-r$, giving $t_0\ge r+2(3-r)\ge3$.

It remains to find four degree-five vertices with at most one boundary
neighbor. Their number is at least $t-r\ge b+2\varepsilon-r$.
Each boundary vertex has at least two neighbors in $X$, so $b\ge6$.
The bound is at least four unless $b=6$, $r=3$ and $\varepsilon=0$.
In that case $P$ is a triangulation, with all three boundary vertices
of degree four. Write $S=\{a,b,c\}$ and let $abx,bcy,caz$ be the three
adjacent faces. The vertices $x,y,z$ are distinct. The cyclic orders
at $a,b,c$ force $xz,xy,yz$, respectively. These six vertices form
an octahedron. Any further vertex lies inside one of its triangular
faces and yields a three-vertex cut. Hence $P$ is the octahedron,
contrary to $d_P(x)\ge5$. We may now choose three vertices counted by
$t_0$ and a fourth degree-five vertex with at most one boundary neighbor.
\end{proof}

Fix $X,S,P$ and $U=\{u_1,u_2,u_3,u_4\}$ as in \Cref{thm:structure}
for a smallest counterexample $G$. Put $A=U\cup N_G(U)$. Since the
only edges added to $P$ have both ends in $S$, \eqref{eq:roots} gives
\begin{equation}\label{eq:local}
 |A\cap S|\le1,\qquad G[A]=P[A].
\end{equation}
We shall also use the fact that every $K_4$ subgraph of $G$ avoids $X$.
Indeed, if a $K_4$ meets $X$, all its vertices lie in $X\cup S$,
so it would be a $K_4$ in $P$.

\section{Transferring the planar configurations}\label{sec:transfer}
\subsection{Cores and local patterns}
We use the pattern names of Blanch\'e, Bonamy, and
Bonichon~\cite[Section 4]{BBB}.
A planar graph $J$ is \emph{almost $4$-connected with respect to $U$}
if no three-vertex cut separates two surviving roots, or two surviving
neighbors of a root in the cut.

A \emph{core} $W$ is a $K_4$ or $M_4$ subdivision rooted on $U$, or a
$W_4$ subdivision with rim roots $U$.
Each root has degree three in $W$.
For $u\in U$, its two \emph{remaining neighbors} form the set
\[
 T_W(u)=\{x:ux\in E(J)\setminus E(W)\}.
\]
A two-path decomposition of $W$ will also be called a two-coloring.
Each root is an endpoint exactly once.
Every nonroot has core degree zero, two, or four.
There is at most one nonroot of core degree four.

A pattern treats one or two roots. In a one-root pattern, the remaining
neighbors are $a,b$.
The $D$ and $T1$ patterns use the pairs $\{w,a\}$ and $\{w,b\}$.
The $T2$ patterns use the same pair $\{a,b\}$ at both roots.
The letters $a,b,w$ denote distinct nonroots.
The adjacency conditions are listed below.
\begin{center}\small
\begin{tabularx}{\textwidth}{lY}
\toprule
Pattern&Adjacency and parity conditions\\\midrule
$\pV$&$ab\notin E(J)$.\\
$\pVp$&$ab\in E(W)$.\\
$\pat{N_e},\pat{N_o}$&$ab\in E(J)\setminus E(W)$.
For $N_e$, choose $a$ of even degree. For $N_o$, both $a,b$ have odd degree.\\
$\pU$&The remaining part is $a-u-v-b$, with $uv\notin E(W)$ and $ab\notin E(J)$.\\
$\pat{D_a}$&$uv\in E(W)$; $wa,wb,ab\in E(J)$; $w,a,b\notin V(W)$.\\
$\pat{D_b}$&$uv\in E(W)$; $wa,wb\in E(J)$; $ab\notin E(J)$; $w\notin V(W)$.\\
$\pat{T1_a},\pat{T1_b}$&$wa,wb\in E(J)$.
For $T1_a$, choose $a$ of even degree. For $T1_b$, both $a,b$ have odd degree.\\
$\pat{T2A_a},\pat{T2A_b}$&$ab\in E(J)$.
For $T2A_a$, choose $a$ of odd degree. For $T2A_b$, both have even degree.\\
$\pat{T2NA_a},\pat{T2NA_b}$&$ab\notin E(J)$.
For $T2NA_a$, choose $a$ of even degree. For $T2NA_b$, both have odd degree.\\
\bottomrule
\end{tabularx}\end{center}
These conditions are only part of the pattern definitions.
We use the position and color conditions from~\cite{BBB}.
The conditions needed for restoration are listed explicitly in
\Cref{tab:positions}, after the position check for R3.

A pattern family partitions $U$ into one-root and two-root groups.
Let $V(C)$ consist of the roots and remaining neighbors of a pattern.
The compatibility condition in~\cite{BBB} is as follows.
If one of two patterns is $\pV$, $\pVp$, or $\pU$, they may share
at most one nonroot. Otherwise they are disjoint.
The core must also have a two-coloring that meets all their color
conditions. Configuration R5 uses one further one-vertex overlap.
We check that overlap directly in \Cref{sec:R5}.

We use the configuration existence results
in~\cite[Lemma 4.4, Claim 4.14, and Lemma 4.16]{BBB} in the following form.
\begin{theorem}\label{thm:patterns}
Let $J$ be planar and almost $4$-connected with respect to four
vertices $U$ of degree five.
The constructions in~\cite[Sections 4.5--4.7]{BBB} give a core $W$,
a two-coloring of $W$, and local patterns covering $U$.
The core is a $K_4$, $M_4$, or $W_4$ subdivision.
The patterns satisfy the position and color conditions specified
in those constructions.
\end{theorem}
We use this planar configuration theorem inside the planar part.
The following arguments transfer its operations to $G$.

\subsection{Matching parity}
A boundary vertex may have different degree parity in $P$ and $G$.
For each $s\in S$, add a leaf at $s$ if needed to make these parities
agree. Denote the resulting planar graph by $\widehat P$.
The leaves are used only to choose patterns.
They do not enter the final reduced graph or its decomposition.
The root neighborhoods and all adjacencies in $A$ stay unchanged, and
\begin{equation}\label{eq:parityP}
 d_{\widehat P}(v)\equiv d_G(v)\pmod2\qquad(v\in A).
\end{equation}
The graph $\widehat P$ is almost $4$-connected with respect to $U$.
Indeed, after at most three vertices are deleted, all surviving
vertices of $P$ remain connected. No new leaf is a root neighbor.
Thus the deletion separates neither two surviving roots nor two
surviving neighbors of a deleted root.

Choose the initial subdivision $Q_0$ in $P$ as follows. If $P$ contains
a $K_4$ subdivision rooted on $U$, choose one with the fewest edges.
Its branch paths are chordless. Otherwise use
\cite[Lemma~4.4]{BBB}, based on Yu's theorem~\cite{Yu}, to obtain a
chordless subdivision of $M_4$ rooted on $U$ with the properties below.
For parallel branches, an edge between their ends that constitutes one
branch is allowed to join the ends of the other.

Call a root neighbor \emph{remaining} if its edge to the root lies
outside $Q_0$. In the $M_4$ case, roots with no branch between them
have no common remaining neighbor. An internal vertex on a branch that
is the only one between its ends cannot be a remaining neighbor of
either of the other two roots. Two roots joined by parallel branches
have at most one common remaining neighbor, and any such neighbor lies
inside a parallel branch joining the opposite pair of roots.
These are the conditions of~\cite[Definition~4.3]{BBB}.

A new leaf cannot occur on a branch of a subdivision rooted on $U$.
It is not a root neighbor either. Thus $Q_0$ satisfies the same
conditions in $\widehat P$, and in the $M_4$ case $\widehat P$ still has
no $K_4$ subdivision rooted on $U$, as required in~\cite[Claim~4.19]{BBB}.
Apply the constructions of \Cref{thm:patterns}, beginning with the changes
that give a strong subdivision in the sense of~\cite[Section~4.3]{BBB}.
Carry out all these steps in $\widehat P$ with the position and color
conditions of the chosen configurations.
Let $W$ be the resulting core.

\begin{lemma}\label{lem:locality}
The constructions can be carried out so that
\begin{equation}\label{eq:locality}
 E(W)\subseteq E(Q_0)\cup E(P[A]).
\end{equation}
Every final pattern and every edge added or deleted by a pattern
uses only vertices of $A$.
\end{lemma}
\begin{proof}
The last assertion follows from the pattern definitions.
For the core, keep the fixed set $A=U\cup N_P(U)$ throughout the
construction. Every current remaining neighbor lies in $A$.
So does the first vertex next to a root on a current branch path.
The changes in~\cite{BBB} have the following form.
\begin{center}\small
\begin{tabularx}{\textwidth}{lY}
\toprule
Construction&Edges newly used by the core\\\midrule
CX1--CX4, Section 4.3&Edges from roots to their neighbors, and edges
between these neighbors. The other sections come from old paths.\\
Routing, Definition 4.20&A root edge is replaced through another
neighbor of the same root.\\
D1--D4, Section 4.5&Root--neighbor edges and routing edges, joined to
sections of old branch paths.\\
J1--J6, Section 4.6&No change, or edges within root neighborhoods,
joined to old path sections.\\
R1--R9, Section 4.7&No change, or root--neighbor edges and two-edge
paths through common neighbors, joined to old path sections.\\
\bottomrule
\end{tabularx}\end{center}
More precisely, D1--D4 add edges $u_iv_i$ with $v_i\in N_P(u_i)$,
and routing edges. J1 and J5 only choose colors.
J2 and J6 add edges from roots to contact vertices.
J3 and J4 may also use an edge between two neighbors of one root.
R1, R2, R5, R6, and R8 do not change the core at this stage.
R3 and R9 cut and rejoin old paths through root--neighbor edges;
a common contact may become the single degree-four nonroot.
R4 and R7 replace a branch by a two-edge path through a common neighbor.

For example, CX1 and CX4 insert a two-edge path through a common
neighbor. CX2 inserts $u_1v_1vu_2$. CX3 inserts the initial section
$u_1vv_2$ or $u_1v_1vv_2$ and keeps the rest of an old branch.
Here $v,v_1,v_2$ are remaining neighbors of the roots.
Routing inserts $uv_2$ and $v_2v_1$, where both $v_1,v_2$ neighbor $u$.
Thus all these edges have both ends in the fixed set $A$.
The D, J, and R changes use the same types of edges, as listed above.
Every other section is taken from a previous core.
Induction over these changes gives \eqref{eq:locality}.
The strong-subdivision step uses CX1--CX4; it adds no further type
of edge.
\end{proof}

\subsection{Replacing a boundary edge}\label{sec:boundary}
We first show that $Q_0[S]$ contains at most one edge. All boundary
vertices are nonroots and have degree at most two in $Q_0$. If there
were exactly two edges, they would form a path of length two, since
$|S|\le3$. Its middle vertex has degree two in $Q_0$, so the path
lies on one branch. Its ends are adjacent in the clique $P[S]$, giving
a chord. The exception for parallel branches concerns only their
roots. If $Q_0[S]$ were a triangle, it would be a component of $Q_0$
containing no root. This proves the claim.

By \eqref{eq:local}, $P[A]$ contains no edge with both ends in $S$.
Hence \Cref{lem:locality} gives the same one-edge bound for the final
core $W$. All pattern operations take place in $A$, so they leave this
possible boundary edge unmarked and unchanged.
If it is an added edge $ab\notin E(G)$, replace it by a shortest
outside path $L$.
For $S=\{a,b,c\}$, choose $L$ in $G-(X\cup\{c\})$.
For $S=\{a,b\}$, choose it in $G-X$.
Its interior avoids $X\cup S$.

This replacement only subdivides one core edge.
It removes no old vertex and changes no edge incident with a root.
The new internal vertices cannot be adjacent to roots in $X$,
so they lie outside $A$.
The induced-graph identity \eqref{eq:local} and the parity agreement in \eqref{eq:parityP} preserve the adjacency and parity conditions of every pattern. A pattern vertex that avoided a core path still avoids
it after replacement. Hence the two-coloring and all position
conditions transfer to $G$.
Continue to write $W$ for this core.
Its only possible part outside $P$ is $L$.

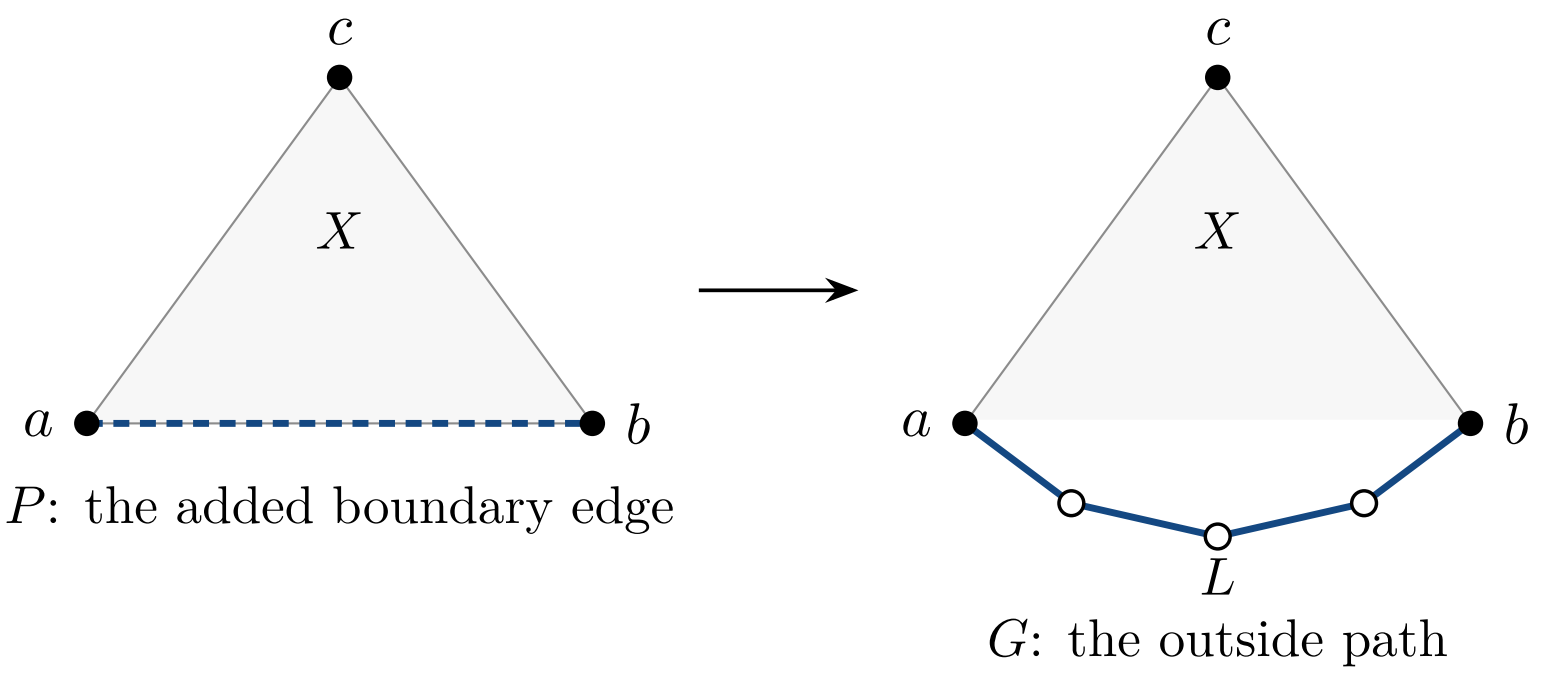
\begin{figure}[tbp]
\centering
\begin{tikzpicture}[scale=.88]
\begin{scope}
 \draw[fill=black!3,draw=black!45](-1.9,0)--(1.9,0)--(0,2.6)--cycle;
 \node[lab]at(0,1.45){$X$};
 \coordinate(a)at(-1.9,0);\coordinate(b)at(1.9,0);\coordinate(c)at(0,2.6);
 \draw[oldpath,densely dashed](a)--(b);
 \foreach \x/\pos in {a/left,b/right,c/above}\node[vertex,label=\pos:$\x$]at(\x){};
 \node[lab]at(0,-.65){$P$: the added boundary edge};
\end{scope}
\draw[arrow](2.7,1)--(3.9,1);
\begin{scope}[xshift=6.6cm]
 \draw[fill=black!3,draw=black!45](-1.9,0)--(1.9,0)--(0,2.6)--cycle;
 \node[lab]at(0,1.45){$X$};
 \coordinate(a)at(-1.9,0);\coordinate(b)at(1.9,0);\coordinate(c)at(0,2.6);
 \draw[white,line width=1.4pt](a)--(b);
 \draw[oldpath](a)--(-1.1,-.6)--(0,-.85)--(1.1,-.6)--(b);
 \foreach \p in {(-1.1,-.6),(0,-.85),(1.1,-.6)}\node[newvertex]at\p{};
 \foreach \x/\pos in {a/left,b/right,c/above}\node[vertex,label=\pos:$\x$]at(\x){};
 \node[lab]at(0,-1.15){$L$};
 \node[lab]at(0,-1.65){$G$: the outside path};
\end{scope}
\end{tikzpicture}
\caption{Replacing a missing boundary edge when $S=\{a,b,c\}$.
The internal vertices of $L$ lie outside $X\cup S$.
The drawing shows the boundary only; the internal core is omitted.}
\label{fig:boundary}
\end{figure}

\subsection{The common neighbors in R3}
We check a position property of the constructions in~\cite{BBB}
before using them for restoration.
\begin{lemma}\label{lem:common}
Every $T1$ or $T2A$ pattern in the final family occurs in R3,
on roots $u_3,u_4$. The other two roots have $V$ patterns.
For $T1$, the common remaining neighbor $w$ lies outside the final
core $W$. For $T2A$, both adjacent common remaining neighbors lie
outside $W$. In particular, there is at most one $\pat{T1_b}$ pattern.
\end{lemma}
\begin{proof}
We first check which configurations use $T1$ or $T2A$.
D1--D4 use one-root patterns.
The two-root patterns in J1--J6 are $U$ or $T2NA$.
R1, R6, and R8 use $D_a,D_b,T2NA$, or one-root patterns.
R2 uses $T2NA$; R4 and R7 use one-root patterns or $U$;
R5 uses $D_b$; and R9 uses one-root patterns or $T2NA$.
Thus only R3 uses $T1$ or $T2A$.
It acts on $u_3,u_4$, while $u_1,u_2$ each have a $V$ pattern.

Let $Q$ be the strong $K_4$ subdivision before the R3 change.
The three core edges at each of $u_3,u_4$ do not change, so their
remaining neighbors stay the same.
Each of these roots has adjacent remaining neighbors and shares one
with the other root. By the definitions in~\cite{BBB}, each root is
unsettled or forms a $V'$ pattern.
The other settled types require nonadjacent remaining neighbors or
no shared remaining neighbor.
Property A in~\cite[Definition 4.11 and Section 4.3]{BBB} rules out
the relevant chords on a branch incident with either root.
Hence a remaining neighbor cannot lie internally on such a branch.
If $w$ lies on $Q$, it must lie inside the $u_1$--$u_2$ branch.

R3 removes that whole branch.
The new core uses other old sections and edges from $u_1,u_2$ to
contact vertices on the old $u_3$--$u_4$ branch.
No internal vertex of the removed branch enters the new core.
If $w$ was outside $Q$, the change does not add it either.
Therefore $w\notin V(W)$.

The same argument applies to each of the two adjacent common neighbors
in a $T2A$ pattern in R3.
A common neighbor on $Q$ must lie on the removed $u_1$--$u_2$ branch;
a common neighbor outside $Q$ is not added by the change.
Thus both lie outside $W$.
The outside replacement path from \Cref{sec:boundary} adds no vertex
of $A$, so these conclusions also hold in $G$.
\end{proof}
Nonadjacent common neighbors may lie on the final core.
In that case we keep the color-avoidance conditions of~\cite{BBB}.

\begin{figure}[tbp]
\centering
\begin{tikzpicture}[scale=.9]
\begin{scope}
 \coordinate(u1)at(-1.65,0);\coordinate(u3)at(0,1.65);
 \coordinate(u2)at(1.65,0);\coordinate(u4)at(0,-1.65);
 \coordinate(h)at(0,0);\coordinate(w)at(0,-2.75);
 \draw[edge](u1)--(u3)--(u2)--(u4)--cycle;
 \draw[edge](u3)--(h)--(u4);
 \draw[newpath](u1)to[out=-100,in=180](w)to[out=0,in=-80](u2);
 \draw[oldpath](u1)--(h)--(u2);
 \draw[pale](w)--(u4);
 \draw[pale](w)..controls(-3,-2.75)and(-3,2.1)..(u3);
 \node[lab,text=pathred]at(1.65,-2.45){$Q_{12}$};
 \foreach \x/\lab/\pos in {u1/u_1/left,u3/u_3/above,u2/u_2/right,u4/u_4/right,h/h/above left,w/w/below}
  \node[vertex,label=\pos:$\lab$]at(\x){};
 \node[lab]at(0,-3.45){before R3};
\end{scope}
\draw[arrow](2.55,-.4)--(3.65,-.4);
\begin{scope}[xshift=6.2cm]
 \coordinate(u1)at(-1.65,0);\coordinate(u3)at(0,1.65);
 \coordinate(u2)at(1.65,0);\coordinate(u4)at(0,-1.65);
 \coordinate(h)at(0,0);\coordinate(w)at(0,-2.75);
 \draw[newpath](u2)--(u4)--(u1)--(h)--(u3);
 \draw[oldpath](u1)--(u3)--(u2)--(h)--(u4);
 \draw[pale](w)--(u4);
 \draw[pale](w)..controls(-3,-2.75)and(-3,2.1)..(u3);
 \foreach \x/\lab/\pos in {u1/u_1/left,u3/u_3/above,u2/u_2/right,u4/u_4/right,h/h/above left,w/w/below}
  \node[vertex,label=\pos:$\lab$]at(\x){};
 \node[lab]at(0,-3.45){the wheel core $W$};
\end{scope}
\end{tikzpicture}
\caption{The wheel case of R3. On the left, remove the dashed branch
$Q_{12}$ from the core and use the two blue edges at $h$.
The right drawing shows a two-path decomposition of the new core.
Here the common remaining neighbor $w$ was on $Q_{12}$, so it is
outside $W$. The pale edges $wu_3,wu_4$ are not core edges.
Curves may represent longer branch paths; crossings are not vertices.
Other edges are omitted.}
\label{fig:R3}
\end{figure}
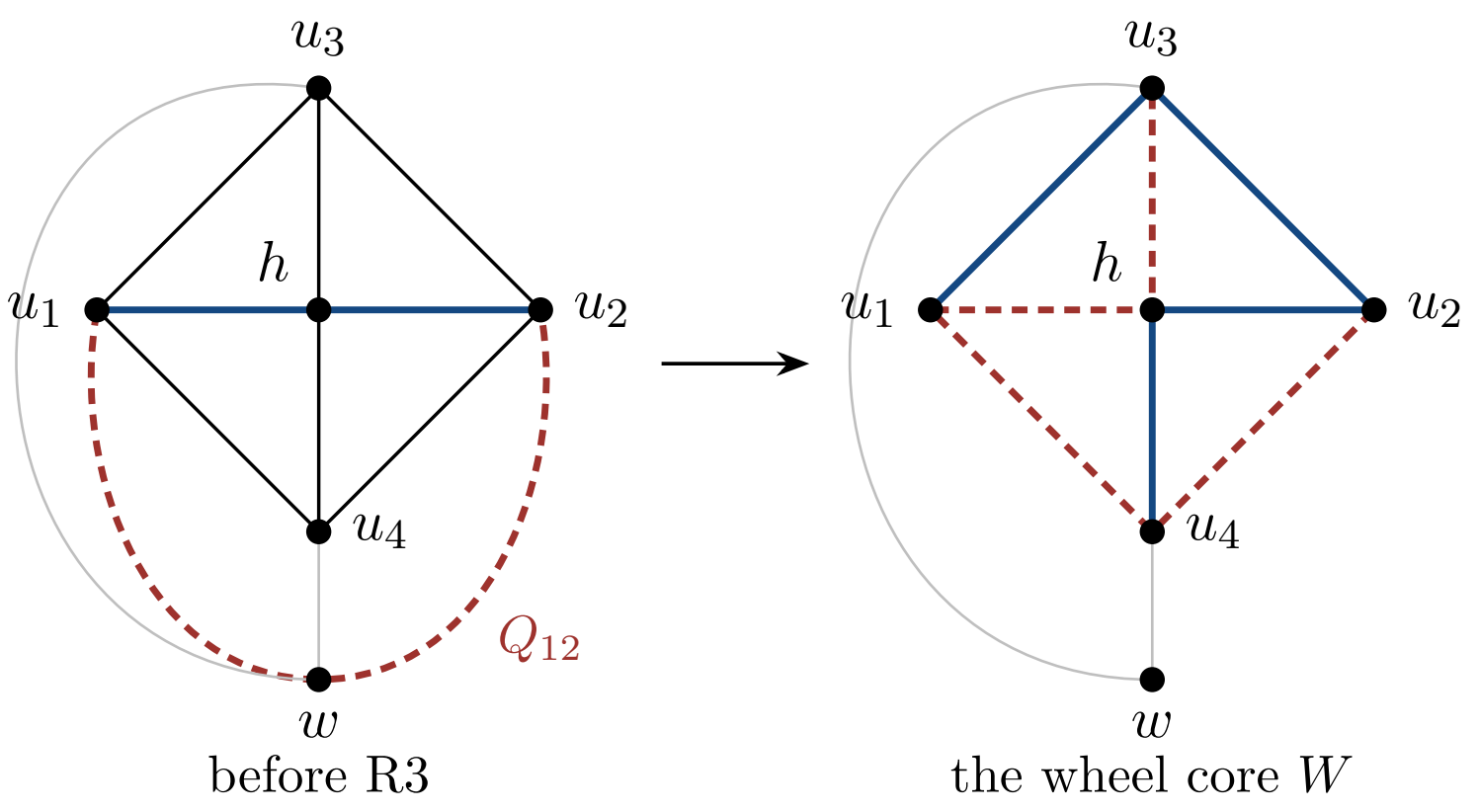

In the wheel case of R3, let $Q_{ij}$ be the old branch between
$u_i$ and $u_j$. The roots $u_1,u_2$ have a common contact $h$ on
$Q_{34}$. The change has the simple form
\[
 E(W)=\bigl(E(Q)\setminus E(Q_{12})\bigr)\cup\{u_1h,u_2h\}.
\]
The old $34$ branch becomes two spokes. The other four old branches
form the rim in the order $u_1,u_3,u_2,u_4$.
Every internal vertex of $Q_{12}$ leaves the core; see \Cref{fig:R3}.

\subsection{Position and color conditions}\label{sec:positions}
Fix the two core paths $P_1,P_2$ before restoration.
For each root $u$, let $P_u$ be the unique one that ends at $u$.
For a core edge $uv$, let $P_{uv}$ be the one containing that edge.
The statement $P_u\ne P_v$ means that the two roots are endpoints
of different core paths.

\Cref{tab:positions} uses the same labels as the adjacency table.
The rows for $T1$ and $T2A$ include the outside-core facts proved in
\Cref{lem:common}. The other requirements are the position and color
conditions of~\cite[Section 4.2]{BBB}, written as vertex exclusions.

\begin{table}[htbp]
\centering\small
\caption{Position and color conditions used in restoration.}
\label{tab:positions}
\renewcommand{\arraystretch}{1.22}
\begin{tabularx}{\textwidth}{lYY}
\toprule
Pattern&Position and color conditions&Change to a core path\\\midrule
$\pV,\pVp,\pU$&No further condition&None\\
$\pat{N_e}$&$b\notin V(P_u)$&Append $ub$ at $u$\\
$\pat{N_o}$&$\{a,b\}\cap V(P_u)=\varnothing$&Append $uba$ at $u$\\
$\pat{D_a}$&$\{w,a,b\}\cap V(W)=\varnothing$&In $P_{uv}$, replace $uv$ by $uwabv$\\
$\pat{D_b}$&$w\notin V(W)$&In $P_{uv}$, replace $uv$ by $uwv$\\
$\pat{T1_a}$&$w\notin V(W)$; $b\notin V(P_u)$&Append $uwb$ at $u$\\
$\pat{T1_b}$&$P_u\ne P_v$; $w\notin V(W)$;
$a\notin V(P_u)$; $b\notin V(P_v)$&Append $uwa$ at $u$ and $vwb$ at $v$\\
$\pat{T2A_a}$&$\{a,b\}\cap V(W)=\varnothing$&Append $uba$ at $u$\\
$\pat{T2A_b}$&$P_u\ne P_v$; $\{a,b\}\cap V(W)=\varnothing$&Append $uba$ at $u$ and $va$ at $v$\\
$\pat{T2NA_a}$&$b\notin V(P_u)$&Append $ub$ at $u$\\
$\pat{T2NA_b}$&$b\notin V(P_u)$; $a\notin V(P_v)$;
$P_u=P_v$ is allowed&Append $ub$ at $u$ and $va$ at $v$\\
\bottomrule
\end{tabularx}
\end{table}
When two tails are added to the same core path, their new vertices
must be distinct. Tails on different paths may meet: for example,
the two tails of $T1_b$ both use $w$.
Different patterns that add tails or replacement sections have
disjoint nonroot sets by compatibility, apart from R5.
That overlap is checked in \Cref{sec:R5}.
The patterns $V,V',U$ add no core section, so their permitted overlaps
do not affect this condition.
\FloatBarrier

\section{The reduced graph and restoration}\label{sec:restore}
Apply the chosen local reductions in $G$.
Delete the four roots and the core edges, and make the neighbor-edge
changes prescribed by the patterns. Let the resulting graph be $H$.
Keep isolated vertices for now. Thus
\begin{equation}\label{eq:orderH}
 |V(H)|=n-4.
\end{equation}
An edge is \emph{marked} if it will be replaced by a path through deleted
roots during restoration. It may be an auxiliary edge or a retained
original edge. Different marked paths have disjoint internal root sets.
Each root enters at most one old member, either inside a marked path
or at an endpoint of an old path.

\begin{lemma}\label{lem:minor}
The graph $H$ is a simple minor of $G$.
\end{lemma}
\begin{proof}
First consider an edge $uv$ between roots that is outside the core.
If $u,v$ belonged to different patterns, those patterns would both
contain $u,v$, contrary to compatibility.
Thus such an edge belongs to one $\pU$ pattern.
We do not need the roots to be independent outside the core.

All new edges arise as follows:
\begin{center}\small
\begin{tabular}{lll}
\toprule
Pattern&New edge&Path in $G$\\\midrule
$\pV,\pVp$&$ab$&$aub$\\
$\pU$&$ab$&$auvb$\\
$\pat{D_b}$&$ab$&$auvb$, with $uv$ a core edge\\
$\pat{T2NA_a}$&$ab$&$avb$, using one root\\
\bottomrule
\end{tabular}\end{center}
Every other pattern only deletes edges or marks a retained edge.
Retain these replacement paths in $G$.
Delete the other required edges, and then suppress the paths.
Their internal roots belong to different patterns, so they are disjoint.
For $\pat{D_b}$, retain $uv$ until this suppression is complete.

The pattern conditions rule out loops and repeated edges.
For $\pVp$, the original edge $ab$ has been deleted with the core before
its auxiliary copy is added. The extra overlap in R5 shares only one
endpoint and creates no repeated edge.
These deletions and contractions produce exactly $H$.
\end{proof}

\subsection{Inserting new vertices}
The planar proof in~\cite[Lemma 4.15]{BBB} already uses mixed
decompositions. We record the insertion step needed for the reductions
in $G$.
Let $\mathcal E$ partition $E(H)$ into simple paths and cycles.
We allow two operations on its members:
replace a marked edge $ab$ by a path $aRb$ with interior in $U$, or
extend an old path ending at $a$ by $au$, where $u\in U$.
No endpoint extension is applied to a cycle.
All root sets used by these operations are pairwise disjoint.

\begin{lemma}\label{lem:insert}
Suppose that the old vertices chosen for endpoint extensions are
distinct and have odd degree in $H$.
Then all the operations above can be made at once.
Old paths remain simple paths, and old cycles remain simple cycles.
The number of members does not change.
Members from distinct components of $H$ remain vertex-disjoint.
\end{lemma}
\begin{proof}
For an old vertex $a$, let $p_a$ count the paths ending at $a$,
$i_a$ the paths through $a$ internally, and $c_a$ the cycles through $a$.
Then
\[
 d_H(a)=p_a+2i_a+2c_a.
\]
If $d_H(a)$ is odd, then $p_a\ge1$.
Thus each required endpoint path exists.

Expand all marked edges first. The inserted vertices were absent
from $H$, and different expansions insert disjoint sets.
No member gains a repeated vertex, and its type is unchanged.
Next make the endpoint extensions. The added roots are new to every
old member. If a path is selected twice, the two choices concern its
two distinct endpoints and insert distinct roots.
So it remains simple. No operation changes the number of members.
Distinct components started with disjoint vertex sets and receive
disjoint sets of new roots. This proves the last assertion.
\end{proof}

\subsection{The local operations}\label{sec:operations}
In the following table, $u,v$ are the roots of the pattern.
The one-root patterns have remaining pair $\{a,b\}$.
The $D,T1$ patterns have pairs $\{w,a\},\{w,b\}$.
The $T2$ patterns have pair $\{a,b\}$ at both roots.
For $\pU$, the remaining part is $a-u-v-b$.
All changes to $H$ are made after deleting $U$ and $E(W)$.
The notation $a\to u$ means to extend an old path at $a$ by $au$.
The last column only lists changes to the fixed two-coloring of the
core. Unlisted sections keep their colors.

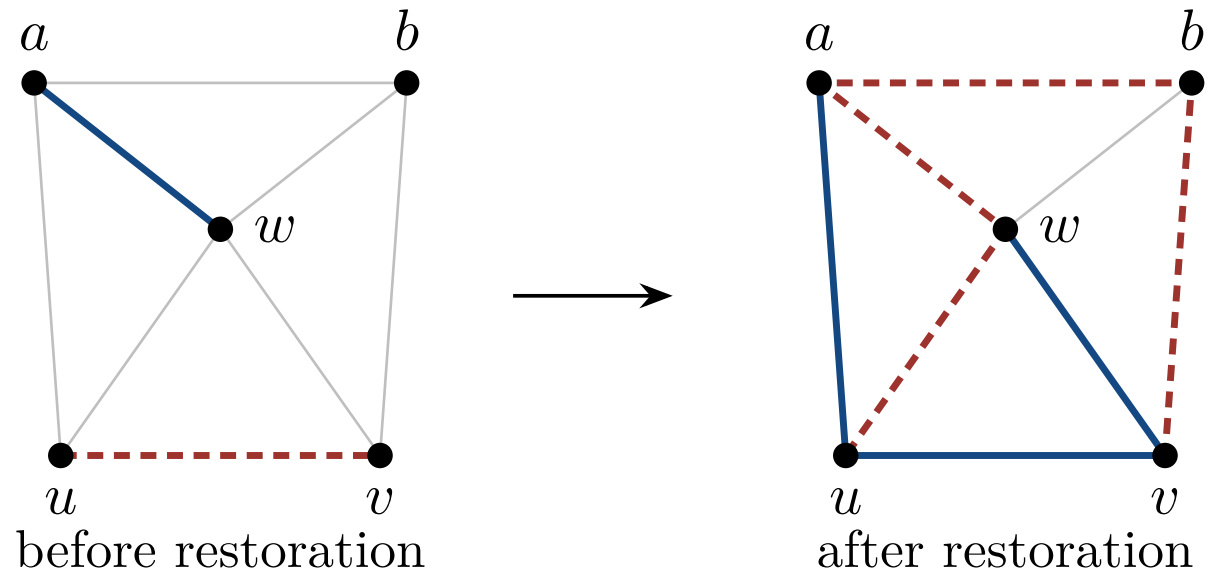
\begin{figure}[tbp]
\centering
\begin{tikzpicture}[scale=.88]
\begin{scope}
 \coordinate(u)at(-1.2,0);\coordinate(v)at(1.2,0);
 \coordinate(w)at(0,1.7);\coordinate(a)at(-1.4,2.8);\coordinate(b)at(1.4,2.8);
 \draw[oldpath](w)--(a);
 \draw[newpath](u)--(v);
 \draw[pale](u)--(w)--(v) (u)--(a)--(b)--(v) (w)--(b);
 \foreach \x/\pos in {u/below,v/below,w/right,a/above,b/above}
  \node[vertex,label=\pos:$\x$]at(\x){};
 \node[lab]at(0,-.7){before restoration};
\end{scope}
\draw[arrow](2.2,1.2)--(3.4,1.2);
\begin{scope}[xshift=5.9cm]
 \coordinate(u)at(-1.2,0);\coordinate(v)at(1.2,0);
 \coordinate(w)at(0,1.7);\coordinate(a)at(-1.4,2.8);\coordinate(b)at(1.4,2.8);
 \draw[oldpath](w)--(v)--(u)--(a);
 \draw[newpath](u)--(w)--(a)--(b)--(v);
 \draw[pale](w)--(b);
 \foreach \x/\pos in {u/below,v/below,w/right,a/above,b/above}
  \node[vertex,label=\pos:$\x$]at(\x){};
 \node[lab]at(0,-.7){after restoration};
\end{scope}
\end{tikzpicture}
\caption{The $\pat{D_a}$ operation.
The solid blue section belongs to an old member.
The dashed red section belongs to one of the two new core paths.
The edge $wb$ is unchanged. All other sections are omitted.}
\label{fig:Da}
\end{figure}

\begin{table}[tbp]
\centering\small
\caption{Local restoration. Every row retains the position and color
conditions of its pattern.}
\label{tab:restore}
\setlength{\tabcolsep}{4pt}
\renewcommand{\arraystretch}{1.22}
\begin{tabularx}{\textwidth}{lYYY}
\toprule
Pattern&Extra change in $H$&Old path or cycle&Two new paths\\\midrule
$\pV$&Add $ab$&$ab\mapsto aub$&No change\\
$\pVp$&Add back core edge $ab$&$ab\mapsto aub$&Keep original $ab$ in its core color\\
$\pat{N_e}$&None&$a\to u$&Append $ub$ at $u$\\
$\pat{N_o}$&Delete $ab$&$a\to u$&Append $uba$ at $u$\\
$\pU$&Add $ab$&$ab\mapsto auvb$&No change\\
$\pat{D_a}$&Delete $ab$; keep $wa,wb$&$wa\mapsto wvua$&Replace core edge $uv$ by $uwabv$\\
$\pat{D_b}$&Add $ab$; keep $wa,wb$&$ab\mapsto auvb$&Replace core edge $uv$ by $uwv$\\
$\pat{T1_a}$&None&$wb\mapsto wvb$; $a\to u$&Append $uwb$ at $u$\\
$\pat{T1_b}$&Delete $wa,wb$&$a\to u$; $b\to v$&Append $uwa,vwb$ to distinct paths\\
$\pat{T2A_a}$&None&$ab\mapsto avb$; $a\to u$&Append $uba$ at $u$\\
$\pat{T2A_b}$&Delete $ab$&$a\to u$; $b\to v$&Append $uba,va$\\
$\pat{T2NA_a}$&Add $ab$&$ab\mapsto avb$; $a\to u$&Append $ub$ at $u$\\
$\pat{T2NA_b}$&None&$a\to u$; $b\to v$&Append $ub,va$\\
\bottomrule
\end{tabularx}
\end{table}
The names $\pat{D_a},\pat{D_b}$ follow the definitions
in~\cite[pp. 43--44]{BBB}. The edge $ab$ is present for $D_a$ and absent
for $D_b$.
For example, in $D_a$, the old member replaces $wa$ by $wvua$.
The new core section $uwabv$ then covers the original edge $wa$.
The original core edge $uv$ moves to the old member.
Thus each edge is used once; see \Cref{fig:Da}.

If a one-root pattern has $ab\in E(W)$, we always use $\pVp$.
Delete $ab$ with the core, add it back in $H$, and expand it to $aub$
in its old member. The original $ab$ keeps its core color.
This avoids treating a deleted core edge as a retained neighbor edge
in an ordinary $\pat{N_e}$ or $\pat{N_o}$ pattern.
The change needs no new core coloring. A $\pVp$ pattern adds no tail.
It only relaxes the compatibility and avoidance conditions of the
corresponding $N$ pattern. Thus all ordinary $N$ patterns used below
have their neighbor edge outside the core.

\Cref{tab:restore} also shows that each root enters exactly one old
member. For $T1_a$ and $T2A_a$, expansion and extension use different
roots. For $D_a,D_b,U$, both roots enter one marked section and are
not used again at an old endpoint.

\subsection{Parity in the final reduced graph}\label{sec:parityH}
For a nonroot $z$, let $\delta(z)$ be the degree loss in addition to
the deletion of core edges. Thus
\begin{equation}\label{eq:degree}
 d_H(z)=d_G(z)-d_W(z)-\delta(z).
\end{equation}
Contributions from different patterns are added.
Marking a retained edge does not change its degree contribution.
\begin{center}\small
\renewcommand{\arraystretch}{1.14}
\begin{tabularx}{\textwidth}{lYY}
\toprule
Pattern&Contribution to $\delta$&Required odd endpoints in $H$\\\midrule
$\pV,\pVp,\pU$&Zero at every nonroot&None\\
$\pat{N_e}$&$\delta(a)=\delta(b)=1$&$a$; originally even\\
$\pat{N_o}$&$\delta(a)=\delta(b)=2$&$a$; originally odd\\
$\pat{D_a}$&$\delta(w)=\delta(a)=\delta(b)=2$&None\\
$\pat{D_b}$&$\delta(w)=2$, $\delta(a)=\delta(b)=0$&None\\
$\pat{T1_a}$&$\delta(w)=2$, $\delta(a)=\delta(b)=1$&$a$; originally even\\
$\pat{T1_b}$&$\delta(w)=4$, $\delta(a)=\delta(b)=2$&$a,b$; originally odd\\
$\pat{T2A_a}$&$\delta(a)=\delta(b)=2$&$a$; originally odd\\
$\pat{T2A_b}$&$\delta(a)=\delta(b)=3$&$a,b$; originally even\\
$\pat{T2NA_a}$&$\delta(a)=\delta(b)=1$&$a$; originally even\\
$\pat{T2NA_b}$&$\delta(a)=\delta(b)=2$&$a,b$; originally odd\\
\bottomrule
\end{tabularx}\end{center}
Every nonroot has even degree in $W$.
Compatibility allows another pattern to meet a required odd endpoint
only if that pattern is $\pV,\pVp$, or $\pU$.
Its contribution is zero. Hence the required endpoints are odd after
all reductions, not just after one local step.
For $\pVp$, deleting the two root edges and adding back $ab$ gives
$\delta(a)=\delta(b)=0$ relative to the core deletion.
The actual degrees may fall, but their parity is preserved.

\subsection{The overlap in R5}\label{sec:R5}
Configuration R5 in~\cite{BBB} uses a $\pat{D_b}$ pattern and an $N$
pattern with one common nonroot. We verify this overlap directly.
Let $D_b$ use roots $u,v$ and remaining vertices $w,a,b$.
Let the $N$ pattern use root $p$ and remaining pair $a,d$.
The vertices $w,a,b,d$ are distinct and outside the core.

The $D_b$ operation adds $ab$, marks it as $auvb$, and replaces the
core edge $uv$ by $uwv$.
The $N$ operation retains or deletes the different edge $ad$.
The contribution of $D_b$ at $a$ is zero, so it does not change the
parity needed by $N$.
The old-member operations use distinct roots $u,v,p$.
The new core section uses $w$, while the $N$ tail uses $a,d$.
Thus the new sections do not conflict.
\Cref{lem:insert} still applies if the two old-member operations act
on the same old path.
If the marked edge lies on an old cycle, the endpoint extension uses
an old path instead. Its existence follows from odd degree.

\begin{theorem}\label{thm:mixed}
Choose the patterns and core coloring by the constructions
of~\cite[Sections 4.5--4.7]{BBB}.
Use $\pVp$ whenever the neighbor edge of a one-root pattern lies in the
core. Use the position and color conditions in \Cref{tab:positions},
and allow the overlap in R5 as checked above.
Every mixed decomposition of $H$ with $m$ members extends to a mixed
decomposition of $G$ with $m+2$ members.
Old paths remain paths, old cycles remain cycles, and the two new
members are simple paths.
Each root belongs to at most one restored old member.
\end{theorem}
\begin{proof}
Compatibility gives distinct marked edges and distinct old endpoints
for extension. By \Cref{tab:restore}, each root is used in only one
old-member operation. The degree table and \eqref{eq:degree} make
all required old endpoints odd.
Thus \Cref{lem:insert} applies to every old path and cycle at once.

For the two new paths, use the fixed core coloring and
\Cref{tab:restore,tab:positions}.
Each added tail avoids the rest of its core path.
In $D_a,D_b$, the replacement section has its interior outside $W$.
Thus each individual change preserves a simple path.

We also check changes made at the same time.
The two tails in $T1_b$ or $T2A_b$ go to different core paths.
In $T2NA_b$, they may go to the same path, but their new vertices
are the distinct vertices $a,b$.
Different patterns that change a core path have disjoint nonroot
sets, except in R5. The $V,V',U$ patterns do not change a core path.
In R5, the replacement section uses $w$, while the other tail uses
only $a,d$; these vertices are distinct by \Cref{sec:R5}.
So two changes on the same core path never give it a repeated vertex.
The two new members are therefore simple paths.

Each row of \Cref{tab:restore} covers its edges exactly once.
All auxiliary edges are expanded, all temporarily deleted original
edges are restored, and every root edge is covered.
Other original edges are unchanged.
The resulting mixed decomposition has $m+2$ members.
\end{proof}
The outside-core conclusion of \Cref{lem:common} is essential here.
Both new tails of $T1_b$ pass through $w$. If $w$ were on the core,
at least one of the two paths would repeat it.

\section{Exceptional components}\label{sec:exceptional}
\subsection{Components isomorphic to \texorpdfstring{$K_5-e$}{K5-e}}
\begin{lemma}\label{lem:marked}
Every $K_5-e$ component of $H$ contains a marked edge.
\end{lemma}
\begin{proof}
Suppose that a $K_5-e$ component $K$ has no marked edge. All edges
of $K$ are original, and none was deleted with $W$: every edge added
back after the core deletion is marked. Each vertex of $K$ belongs
to one of its two $K_4$ subgraphs. The observation following
\eqref{eq:local} therefore gives $V(K)\cap X=\varnothing$.
All changes outside $P$ consist of deleting the edges of the possible
outside path $L$. Thus any vertex of $K$ outside $S$ and $L$ keeps
its original degree, which is at most four.

If there is no outside path, this already gives at least
$5-|S|\ge2$ vertices of degree at most four in $G$.
Suppose that $L$ replaces $ab$ with $a,b\in S$, and put
$S'=S\setminus\{a,b\}$. The path $L$ is shortest in $G-(X\cup S')$.
As $|S'|\le1$, the graph $K-S'$ has diameter at most two and lies
in $G-(X\cup S')$. Hence any two vertices of $K$ on $L$ are at
distance at most two along $L$. If $L$ met $K$ in three vertices,
they would be consecutive on $L$. Neither edge between consecutive
vertices of this triple belongs to $K$, since both were deleted with
the core. But $K$ has only one missing edge. Thus $L$ meets at most
two vertices of $K$. Since $S'$ accounts for at most one more, at
least two vertices of $K$ keep their original degrees. In either
case we contradict \Cref{lem:low}.
\end{proof}

\Needspace{12\baselineskip}
\subsection{A possible triangle component}
\begin{proposition}\label{prop:triangle}
If $H$ has a triangle component $K$, then $W$ is a $W_4$ subdivision
and the pattern family contains a $\pat{T1_b}$ pattern from R3.
Let $h$ be the center of $W$ and $w$ the common remaining neighbor
of this pattern. Then
\[
 V(K)=\{h,w,z\},\qquad d_G(h)=d_G(w)=6,\qquad
 d_G(z)\in\{2,4\},\qquad hw\in E(G).
\]
There is at most one triangle component.
\end{proposition}
\begin{proof}
We first determine which pattern vertices can lie in a triangle
component. For $N,T1_a,T2A,T2NA$, the relevant vertices are either
odd in $H$ or connected to an odd vertex by retained or auxiliary
edges. None can belong to a triangle component.
If a one-root neighbor edge lies in the core, we use $V'$ instead of
$N$, as specified in \Cref{sec:operations}.

A $D_a$ pattern retains $wa,wb$ and deletes $ab$.
If one of these three vertices belongs to a triangle component,
all three must belong to it. But $ab$ is absent, a contradiction.
A $D_b$ pattern retains $wa,wb$ and adds $ab$.
If they form a triangle component, then $w\notin V(W)$ and
$\delta(w)=2$, so $d_G(w)=4$.
Also $\delta(a)=\delta(b)=0$.
At most one of $a,b$ has core degree four.
The other therefore has degree at most four in $G$.
This contradicts \Cref{lem:low}.
At the extra overlap in R5, the connection to the odd endpoint of
the $N$ pattern excludes a triangle component as well.

For $T1_b$, the vertices $a,b$ are odd in $H$.
Only the common neighbor $w$ can lie in a triangle component.
Every other possible vertex has zero additional loss, and hence
original degree $2+d_W(z)$. Such a vertex has degree at most four
in $G$, unless it is the unique nonroot of core degree four.
Let $q$ count the $T1_b$ patterns, let $\eta$ count the degree-four
nonroots of $W$, let $\ell=|V_{\le4}(G)|$, and let $t$ count the
triangle components of $H$.
These components are vertex-disjoint, so
\begin{equation}\label{eq:trianglecount}
 3t\le q+\eta+\ell.
\end{equation}
By \Cref{lem:common}, $q\le1$. Also $\eta\le1$ and $\ell\le1$.
If $t\ge1$, equality is forced throughout:
\[
 t=q=\eta=\ell=1.
\]
The three vertices must be $w$, the unique center $h$, and the unique
vertex $z$ of degree at most four.

By \Cref{lem:common}, $d_W(w)=0$.
The degree table gives $\delta(w)=4$, so $d_G(w)=2+4=6$.
The center $h$ is not one of the odd vertices of the $T1_b$ pattern.
The other two patterns are $V$ and contribute zero.
Thus $d_G(h)=2+4=6$.
The vertex $z$ has core degree zero or two and additional loss zero.
Hence $d_G(z)\in\{2,4\}$.

In the wheel case of R3, the edges $hu_1,hu_2$ are original core edges.
Only the $V$ patterns at $u_1,u_2$ can create auxiliary edges.
Neither remaining pair contains $h$.
Thus every edge of $H$ incident with $h$ is original.
In particular, $hw\in E(G)$.
\end{proof}

\begin{figure}[tbp]
\centering
\begin{tikzpicture}[scale=.9]
\begin{scope}
 \coordinate(u1)at(-1.6,0);\coordinate(u3)at(0,1.6);
 \coordinate(u2)at(1.6,0);\coordinate(u4)at(0,-1.6);
 \coordinate(h)at(0,0);\coordinate(w)at(3.1,0);
 \draw[edge](u1)--(u3)--(u2)--(u4)--cycle;
 \draw[edge](h)--(u1) (h)--(u2);
 \draw[pale](h)--(u3) (h)--(u4);
 \draw[oldpath](w)to[bend right=18](u3) (w)to[bend left=18](u4);
 \draw[newpath](h)to[bend right=29](w);
 \foreach \x/\lab/\pos in {u1/u_1/left,u3/u_3/above,u2/u_2/right,u4/u_4/below,h/h/below left,w/w/right}
 \node[vertex,label=\pos:$\lab$]at(\x){};
 \node[lab]at(.7,-2.3){the forbidden edge $hw$};
\end{scope}
\draw[arrow](4.2,0)--(5.1,0);
\begin{scope}[xshift=7.1cm]
 \foreach \y/\l in {1.6/u_1,0/u_2,-1.6/w}{\coordinate(\l)at(-.8,\y);}
 \foreach \y/\l in {1.6/u_3,0/u_4,-1.6/h}{\coordinate(\l)at(.8,\y);}
 \foreach \a in {u_1,u_2,w}\foreach \b in {u_3,u_4,h}\draw[edge](\a)--(\b);
 \foreach \a in {u_1,u_2,w}\node[vertex,label=left:$\a$]at(\a){};
 \foreach \b in {u_3,u_4,h}\node[newvertex,label=right:$\b$]at(\b){};
 \node[lab]at(0,-2.3){the resulting $K_{3,3}$ subdivision};
\end{scope}
\end{tikzpicture}
\caption{The obstruction in \Cref{lem:wheel}.
Curves represent branch paths; crossings are not vertices.
Delete the two pale spokes and retain the dashed edge $hw$.
The two branch classes are $\{u_1,u_2,w\}$ and $\{u_3,u_4,h\}$.}
\label{fig:wheel}
\end{figure}
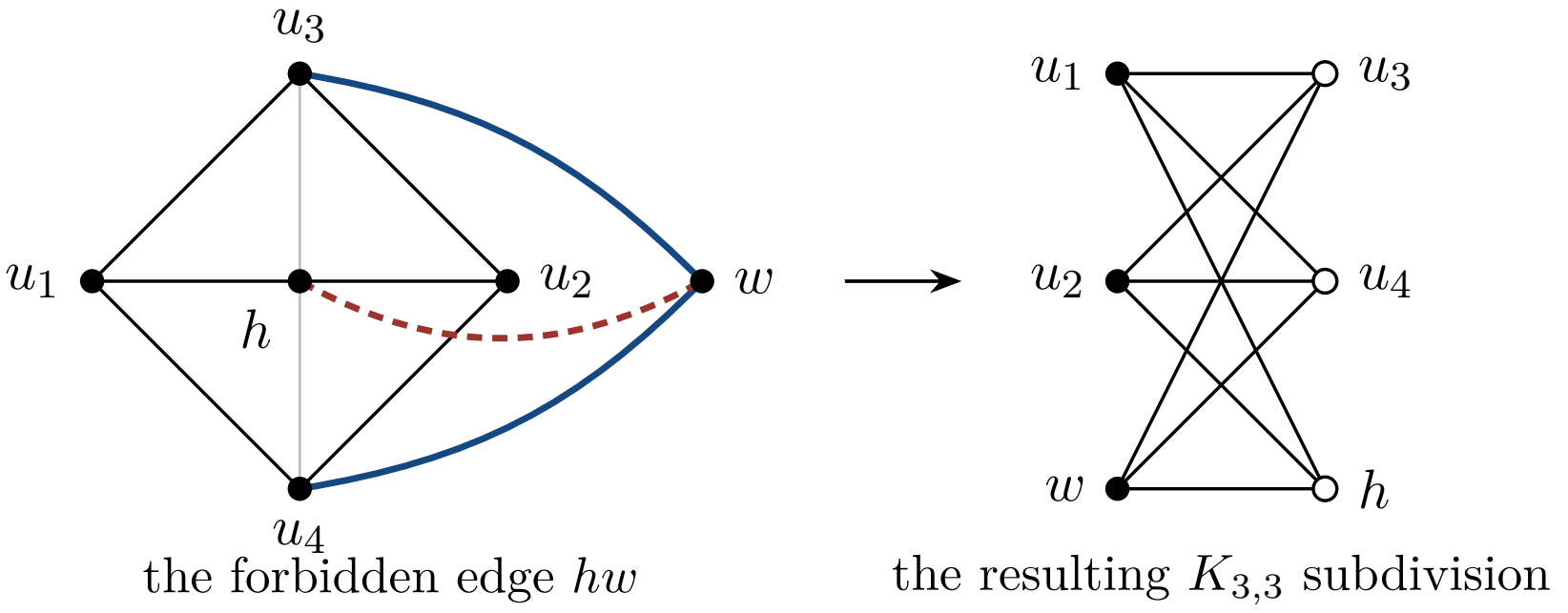

The following observation rules out the wheel in \Cref{prop:triangle}.
\begin{lemma}\label{lem:wheel}
Let $J$ be planar and contain a $W_4$ subdivision with center $h$ and
rim roots in the order $u_1,u_3,u_2,u_4$.
If $w$ is outside the subdivision and $wu_3,wu_4\in E(J)$, then
$hw\notin E(J)$.
\end{lemma}
\begin{proof}
Suppose $hw\in E(J)$. Use $hw$ and the two spokes to $u_1,u_2$ for
the three branches from $h$ to $\{w,u_1,u_2\}$. Each of $u_3,u_4$
has an edge to $w$ and one rim branch to each of $u_1,u_2$.
These nine paths have disjoint interiors, since $w$ is outside the
wheel and the rim branches and spokes belong to one subdivision.
They form a $K_{3,3}$ subdivision in $J$, contrary to planarity
(\Cref{fig:wheel}).
\end{proof}

\begin{theorem}\label{thm:notriangle}
The graph $H$ has no triangle component.
\end{theorem}
\begin{proof}
Suppose that it has one.
By \Cref{prop:triangle}, we are in the wheel case of R3.
Relabel the roots as in that configuration.
The $T1_b$ pattern acts on $u_3,u_4$ and has common neighbor $w$.
The other two roots have $V$ patterns, and
\begin{equation}\label{eq:wheeledges}
 hw,hu_1,hu_2,wu_3,wu_4\in E(G),\qquad w\notin V(W).
\end{equation}
In particular, $hw$ is original, not auxiliary.

Return to the planar part.
Let $W_P$ be the wheel core before a possible added boundary edge
was replaced by $L$.
Replacing $L$ back by that edge preserves the roots, center, and rim
order. The inserted outside vertices are not in $A$, so the local
vertices and edges in \eqref{eq:wheeledges} are unchanged.
The center $h$ sees two roots, and $w$ sees two other roots.
At most one root has a neighbor in $S$ by \Cref{thm:structure}.
Thus neither $h$ nor $w$ lies in $S$.
As neighbors of roots in $X$, they lie in $X\cup S$.
Consequently $h,w\in X$.
All five edges in \eqref{eq:wheeledges} belong to $P$, and
$w\notin V(W_P)$.

The rim order in R3 is $u_1,u_3,u_2,u_4,u_1$.
Its branches are the old $13,23,24,14$ branches.
The center $h$ is the common contact on the old $34$ branch, and
$hu_1,hu_2$ are the newly used original edges.
Apply \Cref{lem:wheel} to $P,W_P,h,w$.
It gives $hw\notin E(P)$, contrary to \eqref{eq:wheeledges}.
\end{proof}

\section{Proof of the main theorem}\label{sec:proof}
\begin{proof}[Proof of \Cref{thm:main}]
Suppose that $G$ is a smallest counterexample. Since $K_1$ and $K_2$
satisfy the bound, \Cref{lem:low,thm:structure} give a planar part
and four roots. Choose the local configurations in \Cref{sec:transfer}
and form the reduced graph $H$ as in \Cref{sec:restore}.
By \Cref{lem:minor,thm:notriangle}, $H$ is a simple $K_5$-minor-free
graph on $n-4$ vertices and has no triangle component.
Every component $J$ other than $K_5-e$ admits a decomposition into at most $\floor{|V(J)|/2}$ paths
by minimality. Isolated vertices contribute nothing.
Decompose each $K_5-e$ component into a Hamilton cycle and a path of
length four. This gives a mixed decomposition with at most
\begin{equation}\label{eq:countH}
 \sum_{J\text{ a component of }H}\floor{|V(J)|/2}
 \le\floor{(n-4)/2}
\end{equation}
members.

Apply \Cref{thm:mixed}. Restore the old members and add two new paths.
The only cycles left come from the $K_5-e$ components of $H$.
Fix such a component $J$.
Let $C,T$ be the restored members from its Hamilton cycle and its
four-edge path.
By \Cref{lem:marked}, $J$ contains a marked edge.
Each root enters at most one old member, so
\[
 V(C)\cap V(T)\subseteq V(J),\qquad |V(C)\cap V(T)|\le5.
\]
If a marked edge lies on the old cycle, then $C$ has length greater
than five. It cannot form the exceptional path--cycle pair.
Otherwise a marked edge lies on the old path and $C$ is still a
$5$-cycle. The expansion removes an edge of $J$ between its five
vertices, and every edge added to $T$ has a new root as an end.
Hence $C\cup T$ has at most eight edges among the vertices of $C$.
Again the exception does not occur.
By \cite[Lemma 2.5]{ZhangLiuHong}, replace $C\cup T$ by two simple paths.

Do this for every $K_5-e$ component.
Members from distinct old components are vertex-disjoint, including
the inserted roots, by \Cref{lem:insert}.
So one replacement does not affect another old cycle.
The member count never increases.
All members are now simple paths, and every original edge is covered
exactly once. Hence
\[
 \pn(G)\le\floor{(n-4)/2}+2=\floor{n/2},
\]
contrary to the choice of $G$.
\end{proof}

\begin{corollary}\label{cor:gallai}
Every connected $K_5$-minor-free graph $G$ satisfies
$\pn(G)\le\ceil{|V(G)|/2}$.
The only connected graphs in this class that fail the floor bound
are $K_3$ and $K_5-e$.
\end{corollary}
\begin{proof}
Use \Cref{thm:main} outside $\Exc$.
For the exceptions, $\pn(K_3)=2$ and $\pn(K_5-e)=3$.
Both meet the corresponding ceiling bound.
\end{proof}

\FloatBarrier

\end{document}